\documentclass[10pt]{amsart}

\usepackage[all]{xypic}
\usepackage{rotating,graphicx}

\usepackage{latexsym}

\usepackage{amssymb}
\usepackage{amsfonts}
\usepackage{amscd}
\usepackage{amsmath,amsthm}
\usepackage{pifont}

 \usepackage{wasysym}
 \usepackage{MnSymbol}

\newtheorem{lemma}{Lemma}[section]

\newtheorem{theorem}[lemma]{Theorem}

\newtheorem{lem}[lemma]{Lemma}
\newtheorem{prop}[lemma]{Proposition}
\newtheorem{thm}[lemma]{Theorem}
\newtheorem{cor}[lemma]{Corollary}

{
}

\theoremstyle{definition}

\theoremstyle{remark}

\numberwithin{equation}{section}

\newenvironment{pf}{\noindent{\bf Proof.}}{\hfill $\square$\medskip}

\def\CC{{\mathbb C}}

\def\PP{{\mathbb P}}

\def\ZZ{{\mathbb Z}}

\def\0ol{{\bar 0}}
\def\1ol{{\bar 1}}
\def\2ol{{\bar 2}}
\def\ol2{{\bar 2}}
\def\3ol{{\bar 3}}
\def\4ol{{\bar 4}}
\def\5ol{{\bar 5}}
\def\6ol{{\bar 6}}
\def\7ol{{\bar 7}}
\def\8ol{{\bar 8}}
\def\9ol{{\bar 9}}

\def\bold0{{\bf 0}}
\def\bold1{{\bf 1}}
\def\bold2{{\bf 2}} 
\def\bold3{{\bf  3}}
\def\bold4{{\bf 4}}
\def\bold5{{\bf 5}}
\def\bold6{{\bf 6}}
\def\bold7{{\bf 7}}
\def\bold8{{\bf 8}}
\def\bold9{{\bf 9}}

\def\P2Skly{\PP^2_{Skly}}

\def\Alt{\operatorname {Alt}}

\def\Ext{\operatorname {Ext}}

\def\GL{\operatorname {GL}}

\def\Hom{\operatorname {Hom}}

\def\pd{{\operatorname {\partial}}}

\def\sgn{\operatorname {sgn}}

\def\th{\operatorname {th}}    % for writing n^{th}

\def\fchar{\operatorname{char}}

\def\det{\operatorname{det}}

\def\dim{\operatorname{dim}}

\def\extd{\operatorname{d\!}}
\def\Ext{\operatorname{Ext}}

\def\Fract{\operatorname{Fract}}

\def\GKdim{\operatorname{GKdim}}
\def\gldim{\operatorname{gldim}}

\def\Hom{\operatorname{Hom}}

\def\id{\operatorname{id}}

\def\pd{{\partial}}

\def\Projnc{\operatorname{Proj}_{nc}}

\def\rank{\operatorname{rank}}

\def\Sym{\operatorname{Sym}}

\def\ul1{\operatorname{\underline{1}}}

\def\l{\leftarrow}

\def\d{\downarrow}

\def\a{\alpha}
\def\b{\beta}
\def\c{\gamma}
\def\d{\delta}

\def\l{\lambda}

\def\s{\sigma}

\def\ve{\varepsilon}

\def\sM{{\sf M}}

\def\sT{{\sf T}}

\def\sfd{{\sf d}}

\def\sfu{{\sf u}}

\def\sfv{{\sf v}}
\def\sfw{{\sf w}}
\def\sfx{{\sf x}}

\def\dirlim{\mathop{\vtop{\baselineskip -100pt\lineskip -1pt\lineskiplimit 0pt
\setbox0\hbox{lim}\copy0\hbox to \wd0{\rightarrowfill}}}\limits}
\def\invlim{\mathop{\vtop{\baselineskip -100pt\lineskip -1pt\lineskiplimit 0pt
\setbox0\hbox{lim}\copy0\hbox to \wd0{\leftarrowfill}}}\limits}

\def\I11{{1 \kern -0.8pt \! \mbox{l}}}
\def\mumu{{\mu\kern-4.2pt\mu}}
\def\bfmu{{\mu\kern-4.2pt\mu}}
\def\2slash{\backslash \! \backslash}

\def\boxtimes{\setbox0\hbox{$\Box$}\copy0\kern-\wd0\hbox{$\times$}}

\begin{document}

\title[3-Calabi-Yau algebras]{The Classification of 3-Calabi-Yau algebras with
\\ 3 generators and  3 quadratic relations }

\author{Izuru Mori and S. Paul Smith}

\address{Department of Mathematics, Graduate School of Science, Shizuoka University, Shizuoka 422-8529, Japan.}
\email{simouri@ipc.shizuoka.ac.jp}
\address{Department of Mathematics, Box 354350, Univ.  Washington, Seattle, WA 98195,USA.}
\email{smith@math.washington.edu}

\keywords{Calabi-Yau algebras, Artin-Schelter regular, graded rings, quadratic algebras, Poisson brackets.}

\subjclass{16E65, 16S37, 16S38, 16S80, 16W50}

\thanks{The first author was supported by a Japanese Grant-in-Aid for Scientific Research (C) 91540020. }

\begin{abstract}

Let $k$ be an algebraically closed field of characteristic not 2 or 3, $V$ a 3-dimensional vector space over $k$, $R$ 
a 3-dimensional subspace of $V \otimes V$, and $TV/(R)$ the quotient of the tensor algebra on $V$ by the ideal generated by $R$.
Raf  Bocklandt proved that if $TV/(R)$ is 3-Calabi-Yau, then it is isomorphic to $J(\sfw)$, the ``Jacobian algebra'' of some $\sfw \in V^{\otimes 3}$. This paper classifies the $\sfw\in V^{\otimes 3}$ such that $J(\sfw)$ is 3-Calabi-Yau. 
The classification depends on how $\sfw$ transforms under the action of the symmetric group $S_3$ on $V^{\otimes 3}$ and on  the nature of the subscheme $\{\overline{\sfw}=0\} \subseteq \PP^2$ where $\overline{\sfw}$ denotes the image of $\sfw$ in the symmetric algebra $SV$. 
% Surprisingly, as $\sfw$ ranges over $V^{\otimes 3}-\{0\}$, 
%only nine isomorphism classes of algebras appear as non-3-Calabi-Yau $J(\sfw)$'s. 
%If $J(\sfw)$ is 3-Calabi-Yau, then it is a 3-dimensional Artin-Schelter regular algebra. 
%If $J(\sfw)$ is commutative, it is a polynomial ring on three variables. Thus, the 3-Calabi-Yau algebras of the form $TV/(R)$ are very good
%non-commutative analogues of the polynomial ring. The 3-Calabi-Yau $J(\sfw)$'s tend to arise in families that are flat deformations of the polynomial ring on three variables. 

\end{abstract}

\maketitle
%\tableofcontents

\section{Introduction}

\subsection{}
This paper is motivated by Raf Bocklandt's result  that every graded 3-Calabi-Yau algebra that is the quotient of the path algebra of a quiver is the Jacobian algebra, $J(\sfw)$, of a superpotential, $\sfw$. If $J(\sfw)$ is 3-Calabi-Yau we say that $\sfw$ is {\sf good}. At present, it is not known which  $\sfw$ are good. This paper determines the good $\sfw$ in an important special case, namely that where $J(\sfw)$ is a connected graded algebra on 3 generators subject to
3 homogeneous quadratic relations.

%The ``best'' non-commutative analogues of $\PP^2$ have homogeneous coordinate rings that
%are 3-dimensional Artin-Schelter regular algebras with Hilbert series $(1-t)^{-3}$. 
%If $A$ is such an algebra, the category that plays the role of the category of quasi-coherent sheaves is $\QGr(A)$, the quotient of the category of graded $A$-modules by the Serre subcategory consisting of 
%the graded modules that are the sum of their finite dimensional submodules. 
%There is increasing evidence that for every such  $A$,  there is a 3-Calabi-Yau algebra $A'$ such that $A$ and $A'$ have equivalent categories of graded modules. %$\Gr(A) \equiv \Gr(A')$ and hence $\QGr(A) \equiv \QGr(A')$.
% If that is the case, then understanding {\it all} non-commutative analogues of $\PP^2$ reduces to understandin those arising from the 3-Calabi-Yau algebras of the form studied in this paper.  

\subsection{}
Throughout this paper, $k$ is an algebraically closed field with $\fchar(k) \ne 2,3$, 
$V$ is a 3-dimensional $k$-vector space, $TV$ and $SV$ are the tensor and symmetric algebras on $V$, and 
$R$ is a 3-dimensional subspace of 
$V ^{\otimes 2}$.  We treat $TV$ and $SV$ as graded $k$-algebras with $\deg(V)=1$, and give $TV/(R)$ the induced grading.

\subsection{Calabi-Yau algebras}
Let  $A$ be a graded $k$-algebra, $A^\circ$ its opposite algebra, and $A^e=A \otimes A^\circ$.
We consider $A$ as a left $A^e$-module via $(a \otimes b^\circ)\cdot c=acb$. If $\nu$ is an automorphism of $A$ we write ${{}_{{}_\nu}}\! A_{{}_1}$ for the left $A^e$-module that is $A$ as a graded vector space with action $(a \otimes b^\circ)\cdot c=\nu(a)cb$. 
 We say $A$ is  {\sf twisted Calabi-Yau} of {\sf dimension} $d$ if it has a finite-length resolution as a left $A^e$-module by finitely 
 generated projective $A^e$-modules and there is an isomorphism
$$
\Ext^i_{A^e}(A,A^e) \cong \begin{cases} 
					{{}_{{}_\nu}} A_{{}_1}(\ell) & \text{if $i=d$}
					\\
					0 & \text{if $i \ne d$}
				\end{cases}
$$
of graded right $A^e$-modules for some integer $\ell$ and some graded $k$-algebra automorphism $\nu$. We call $\nu$ the {\sf Nakayama automorphism} of $A$. (Some authors call $\nu^{-1}$ the Nakayama automorphism.) We say $A$ is  {\sf Calabi-Yau} if it is twisted Calabi-Yau and 
$\nu=\id_A$. If $A$ is  twisted Calabi-Yau of dimension $d$, then it has global dimension $d$ by \cite [Prop. 4.5]{YZ}. 

When $A$ is Calabi-Yau of dimension $d$ we often say it is $d$-Calabi-Yau or $d$-CY.

\subsection{Jacobian algebras}
Fix a basis $\{x,y,z\}$ for $V$. The {\sf cyclic partial derivative} with respect to $x$ of a word $\sfw$ in the letters $x,y,z$, 
 is 
$$
 \pd_x (\sfw) \; := \; \sum_{\sfw=\sfu x \sfv} \sfv\sfu
$$
where the sum is taken over all such factorizations. We extend $ \pd_x$ to  $TV$ by linearity.
We define $ \pd_y$ and $ \pd_z$ in a similar way.  The  {\sf Jacobian algebra} associated to $\sfw \in TV$ is 
$$
J(\sfw)  \; :=\; \frac{TV}{( \pd_x (\sfw), \,  \pd_y (\sfw) , \,  \pd_z  (\sfw) )}   \; ,
$$
i.e., $TV$ modulo the ideal generated by the cyclic partial derivatives.
The linear span,
\begin{equation}
\label{defn.Rw}
R_\sfw \; :=\; {\sf span}\{\pd_x(\sfw),\, \pd_y(\sfw),\,\pd_z(\sfw)\},
\end{equation}
does not depend on the choice of basis for $V$ (Lemma \ref{lem.pd.cpd}).

\subsubsection{}
The following is a special case of a result due to Bocklandt.

\begin{thm}
\cite[Thm. 3.1]{B}
\label{thm.raf}
Every connected graded quadratic 3-Calabi-Yau algebra on three degree-one generators 
 is isomorphic to $J(\sfw)$ for some $\sfw \in V^{\otimes 3}$. 
\end{thm}

The problem of classifying those $\sfw \in V^{\otimes 3}$ for which $J(\sfw)$ is 3-Calabi-Yau has been 
of interest since Bocklandt's result (e.g., \cite {BS}).  
This paper solves the problem.

\subsection{The classification}
\label{ssect.main.result}

Suppose $\dim_k(V)=3$ and let $\sfw \in V^{\otimes 3}-\{0\}$. We write $\overline{\sfw}$ for the image of $\sfw$ in $S^3V$, and
$E \subseteq \PP^2$ for the scheme-theoretic zero locus of $\overline{\sfw}$. If $\overline{\sfw}\ne 0$, $E$ is a cubic divisor.
The classification of those $\sfw$ for which $J(\sfw)$ is 3-Calabi-Yau is given in Table \ref{main.table} which we explain 
in \S\S\ref{ssect.hessian}--\ref{ssect.E'} below.

We use the  symbols $\vert\!\vert\!\vert\,$, $\leftrightline \!\!\!\!\!\!\vert\!\vert\,\,$, $\varhexstar\,$, 
$ \vert\!  \largecircle$, $\bigcurlywedge$, $\largetriangleup$, $\varnothing\,$, and $\propto$, 
to denote the singular cubic divisors in $\PP^2$.
%%%  commonly denoted by these symbols.
 
\begin{table}[htdp]
\begin{center}
\begin{tabular}{|l||l|l|}
\hline
$\phantom{xxxxxxx}E$ &  \S\ref{ssect.cw=sw}  \quad $c(\sfw)=s(\sfw)\quad $  & \S\ref{sect.cw.not-=.sw}  $\phantom{xxxx}c(\sfw)\ne s(\sfw) \phantom{\Big)}$   \phantom{xx} 
\\
\hline
\hline
$\PP^2$ or $\vert\!\vert\!\vert\,$, $\leftrightline \!\!\!\!\!\!\vert\!\vert\,\,$, $\varhexstar\,$, 
$ \vert\!  \largecircle$,  $\bigcurlywedge$   & never CY   & always CY $\phantom{\Big)^3}$
\\
\hline
 & &  $E=\largetriangleup.\,$ CY $\Longleftrightarrow H(\overline{\sfw})\ne 8\mu(\sfw)^2\overline{\sfw}  \phantom{\Big)^3}$ \\ 
 $\largetriangleup$, or $\varnothing$, or $\propto$  &   always CY  &  $E=\varnothing. \phantom{x}$ CY $\Longleftrightarrow E'$ is not a triangle  \\
 &  &  $E= \propto.\phantom{x}$ CY $\Longleftrightarrow E'$ is not a triangle  \\
&  &  $\phantom{xxxx} E'$ is defined in \S\ref{ssect.E'}   \\
\hline
 a smooth cubic & CY $\Leftrightarrow j(E) \ne 0$  & always CY   $\phantom{\bigg)}$  \\
\hline
\end{tabular}
\end{center}
\vskip .12in
\caption{When is $J(\sfw)$ 3-Calabi-Yau?}
\label{main.table}
\end{table}

\subsubsection{Hessians}
\label{ssect.hessian}
We write $\nabla^2(f)$ for the matrix of second partial derivatives of $f(x,y,z)$, and $H(f)$ for  
the determinant of $\nabla^2(f)$. We call $H(f)$ the {\it Hessian} of $f$, and write $H^2(f)$ for the Hessian of $H(f)$. 
Whether or not the Hessian of a function in $SV$ is zero does not depend on the choice of basis for $V$ \cite[p.67]{Fischer}. 

Let $f \in S^3V$. Let $E \subseteq\PP^2$ be the  scheme-theoretic zero locus of $f$. 
 
The following facts are well-known and easy to verify.
(Recall that $k=\overline{k}$ and $\fchar(k) \ne 2,3$.)

\begin{itemize}
\item{}
$H^2(f)=0$ if and only if $E$ is $\PP^2$, or a triple line, or a double line and a line, or three lines having a common 
intersection point, or a smooth conic and a line tangent to it, or a cuspidal cubic. 
\item{}
$H^2(f) \ne 0$ if and only if $E$ is a smooth cubic, or three lines with no common point, 
or a smooth conic and line meeting at two points, or a nodal cubic.  
\end{itemize}

\subsubsection{}
The symmetric group $S_3$ acts on $V^{\otimes 3}$ by permuting the components. 
The $c$ and $s$ in Table \ref{main.table} are elements in the group algebra $kS_3$, namely
$$
s\,:=\, \hbox{$\frac{1}{6}$} \sum_{g \in S_3} g
\qquad \hbox{and} \qquad
c \,:=  \,\hbox{$\frac{1}{3}$}\big(1\, + \, (123) \, + \, (321)\big).
$$
We usually think of $c$ and $s$ as linear maps $V^{\otimes 3} \to V^{\otimes 3}$.

\subsubsection{The number $\mu(\sfw)$}
\label{ssect.mu}
Let $\{x,y,z\}$ be a basis for $V$.
The sign representation for $S_3$ occurs in $V^{\otimes 3}$ with multiplicity 1. The element
$$
\sfw_0:=\hbox{$\frac{1}{3}$}\big(xyz+yzx+zxy-zyx-xzy-yxz\big)
$$
is a basis for it. There is a linear map $\mu:V^{\otimes 3} \to k$ such that $c(\sfw)-s(\sfw)=\mu(\sfw)\sfw_0$ for all 
$\sfw \in V^{\otimes 3}$ (see \S\ref{ssect.w0}).

\subsubsection{}
\label{ssect.cancel}
Although $H(\overline{\sfw})$, $\mu(\sfw)$, and $\sfw_0$, depend on the choice of basis for $V$ those dependencies 
 ``cancel out'' in all relevant situations: if $\theta \in \GL(V)$ and 
$H(\overline{\sfw})'$ and $\mu(\sfw)'$ denote the Hessian and $\mu$ computed with respect to the basis 
$\{\theta(x),\theta(y),\theta(z)\}$, then $H(\overline{\sfw})'=(\det \theta)^2H(\overline{\sfw})$ and $\mu(\sfw)'=(\det \theta)\mu(\sfw)$
so the condition $H(\overline{\sfw})\neq 8\mu(\sfw)^2\overline{\sfw}$ does not depend on the choice of basis for $V$.

\subsubsection{The subscheme $E' \subseteq \PP^2$}
\label{ssect.E'}
If $\overline{\sfw} \ne 0$, let $E' \subseteq \PP^2$ denote the zero locus of $H(\overline{\sfw}) + 24 \mu(\sfw)^2 \overline{\sfw}$.
By the remarks in \S\ref{ssect.cancel}, $E'$ does not depend on the basis for $V$  so the condition in Table \ref{main.table} that $E'$
 is or is not a triangle is a basis-free condition.

Furthermore, for a fixed $\overline{\sfw}$ such that $E$ is $\largetriangleup$, $\varnothing$ or $\propto$, there are exactly two values of $\mu(\sfw)$
for which $J(\sfw)$ is not 3-CY. 
 
When $A$ is a 3-dimensional quadratic Artin-Schelter regular algebra,\footnote{The definition is in \S\ref{ssect.AS-reg}.} its moduli space of point modules, often called its {\it point scheme}, plays a central role. It is the largest closed commutative subscheme of $\Projnc(A)$,
the non-commutative $\PP^2$ with homogeneous coordinate ring $A$. The next result is proved in \S\ref{sect.pt.var}.

\begin{thm}
\label{thm.pt.scheme}
When $J(\sfw)$ is 3-Calabi-Yau it  is a 3-dimensional quadratic Artin-Schelter regular algebra and
its point scheme is $E'$. 
\end{thm}

\subsubsection{Effectiveness of the classification}
The conditions in Table \ref{main.table} are effective. 
Given $\sfw$, it is routine to determine if $H^2(\overline{\sfw})$ is zero, and routine to determine if $E$ is singular.
It is easy to determine whether $c(\sfw)$ and $s(\sfw)$ are equal or not. When $E$ is $\varnothing$ or 
$\propto$ it is easy to decide if $E'$ is a triangle of not: first determine $h:=H(\overline{\sfw}) + 24 \mu(\sfw)^2 \overline{\sfw}$, then
use the fact that the zero locus of $h$ is a triangle if and only if $H(h)$ is a non-zero multiple of $h$ \cite[Prop.4.5, p.68]{Fischer}. 
The other thing one needs to determine is $j(E)$ when $E$ is a smooth curve; i.e., given $f \in S^3V$
(in this paper $f=\overline{\sfw}$) compute $j(E)$; solutions to this problem date back to the 1800's.

\subsubsection{}
If $\sfw=0$, then $J(\sfw)=TV$. $TV$ is not 3-Calabi-Yau. 

\subsection{}
Before discussing the proof we need a little more notation and terminology.

We write $\PP^2$ for $\PP(V^*)$ and  $N^{\sT}$ for the transpose of a matrix $N$. 

\subsubsection{Standard algebras and the matrix $Q$}
\label{ssect.standard}
Let $\{x_1,x_2,x_3\}$ be a basis for $V$ and $\{f_1,f_2,f_3\}$ a basis for $R$. Define $\sfx:=(x_1,x_2,x_3)^\sT$ and ${\sf f}:=(f_1,f_2,f_3)^\sT$. There is a unique $3 \times 3$ matrix $\sM$ with entries in $V$
such that ${\sf f}=\sM\sfx$. 

Following \cite[p.34]{ATV1}, we say $TV/(R)$ is {\sf standard} if there are bases for 
$V$ and $R$ such that the entries in $\sfx^\sT \sM$
are also a basis for $R$. In that case, there is a unique $Q \in \GL(3)$ such that $\sfx^\sT \sM =(Q\sM\sfx)^\sT$.

\subsubsection{The symmetrization map $f \mapsto \widehat{f}$}
We write $\Sym^m(V)$ for the space of 
degree-$m$ symmetric tensors and $S^mV$ for the degree$-m$ component of $SV$. The canonical homomorphism $TV \to SV$,
$\sfw \mapsto \overline{\sfw}$,  
restricts to an isomorphism $\Sym^m(V) \to S^mV$ for all $m$; we denote its inverse by $f \mapsto \widehat{f}$.

\subsection{The method of proof}
The first step in our classification is Bocklandt's Theorem.
The second step is to show that $J(\sfw)$ is 3-Calabi-Yau if and only if it is a
3-dimensional Artin-Schelter regular algebra.
We then use the criterion in Theorem \ref{thm.ATV1.1} to decide when $J(\sfw)$ is a 3-dimensional Artin-Schelter regular algebra.

\subsubsection{Artin-Schelter regular algebras}
\label{ssect.AS-reg}

A graded algebra  $A=A_0 \oplus A_1 \oplus \cdots$ is {\sf connected}  if $A_0=k$. 
In that case we may consider $k=A/A_{\ge 1}$ as a left, and as a right, graded $A$-module concentrated in degree zero. 
We say that $A$ is  {\sf Artin-Schelter regular} of dimension $d$ if it has finite GK-dimension, $\gldim(A)=d<\infty$, and
\begin{equation}
\label{AS-Gor-cond}
\Ext^i_A(k_A,A) \cong \Ext^i_A({}_Ak,A) \cong \begin{cases} 
					k(\ell) & \text{if $i=d$}
					\\
					0 & \text{if $i \ne d$.}
				\end{cases}
\end{equation}
for some $\ell \in \ZZ$ (see \cite {AS}). 
Commutative Artin-Schelter regular algebras are polynomial rings.

\begin{thm}
[Corollary \ref{cor.3ASR=3CY}]
\label{thm.3ASR=3CY}
Let $V$ be a 3-dimensional vector space and  $\sfw \in V^{\otimes 3}$. 
Then $J(\sfw)$ is 3-Calabi-Yau if and only if it is a 3-dimensional Artin-Schelter regular algebra. 
\end{thm}

The fact that an Artin-Schelter regular algebra is twisted Calabi-Yau  is proved in \cite[Lem. 1.2]{RRZ}.  In particular, if $J(\sfw)$ is Artin-Schelter regular, then it is twisted Calabi-Yau.  However, Theorem \ref{thm.3ASR=3CY} says that $J(\sfw)$ is not just twisted Calabi-Yau, but Calabi-Yau.  
Although  \cite[Lem. 1.2]{RRZ} is an ``if and only if'' result we can not appeal to it for the other implication in 
Theorem \ref{thm.3ASR=3CY} because we require Artin-Schelter regular algebras to have finite 
GK-dimension whereas \cite{RRZ} does not. Nevertheless,  Corollary \ref{cor.3ASR=3CY} shows that if one of our $J(\sfw)$ is 3-CY, 
then it has finite GK-dimension ({\it that} is the new content of Theorem \ref{thm.3ASR=3CY}) and is therefore Artin-Schelter regular in our sense. 
%For a Jacobian algebra, the following holds in dimension 3.  
 
After proving Theorem \ref{thm.3ASR=3CY}  in \S\ref{sect.prelim}, our strategy is to determine those $\sfw$ for which $J(\sfw)$ 
is Artin-Schelter regular of dimension 3. The main criterion for that is Theorem \ref{thm.ATV1.1} which is due to Artin, Tate, and Van den Bergh.

 Let $\sM$ and $\sfx$ be as in \S\ref{ssect.standard}.
Let $\sfx'$ be another $3 \times 1$ matrix whose entries are a basis for $V$ and $\sM'$ another $3 \times 3$ matrix such that the 
entries in $\sM'\sfx'$ are a basis for $R$.  If we treat $\sM'$ and $\sM$ as matrices with entries in $SV$, then the $2 \times 2$ minors of $\sM$ and $\sM'$ belong to $S^2V$. The rank of $\sM$ is $\ge 2$ at all points in $\PP^2$ if and only if the same is true of $\sM'$ because $\sM'=P_1\sM P_2$ for some $P_1,P_2 \in \GL(3)$. 

\begin{thm}
\cite[Thm.1]{ATV1}
\label{thm.ATV1.1}
Let $V$ be a 3-dimensional vector space, $R$ a 3-dimensional subspace of $V^{\otimes 2}$. Then
$TV/(R)$ is a 3-dimensional Artin-Schelter regular algebra if and only if it is standard in the sense of \S\ref{ssect.standard} and the 
common zero locus of the $2 \times 2$ minors of the matrix $\sM$ in \S\ref{ssect.standard} 
is empty. 
\end{thm}

 \subsubsection{Quadratic Poisson brackets and the algebras $A_{{}_{\pi,\l}} $ and $A_{\!{}_{f,\l}}$}

A {\sf bivector field} on a smooth manifold $M$ is a section of $\wedge^2(TM)$. 
If $\pi$ is a bivector field with the property that the Nijenhuis-Schouten bracket, $[\pi,\pi]_{{}_{NS}}$, is zero, then
the formula
$$
\{f,g\}_{{}_\pi}  \; := \; \langle \pi,\extd f \wedge \extd g\rangle, \qquad f,g \in C^\infty(M),
$$
defines a Poisson bracket on $C^\infty(M)$. Every Poisson bracket on $C^\infty(M)$ is equal to $\{\cdot,\cdot\}_{{}_\pi}$
for a unique bivector field $\pi$.

A {\sf quadratic} Poisson bracket on $SV$ is one such that $\{V,V\}\subseteq S^2V$. 
Only quadratic Poisson brackets appear in this paper.  

For each quadratic Poisson bracket $\{\cdot,\cdot\}_{{}_\pi}$ on $SV$ and every  $\l \in k$ we define
\begin{align*}
\label{defn.Rlambda}
 R_{{}_{\pi,\l}}    \; :=\; \{u \otimes v &- v \otimes u - \l\widehat{\{u,v\}}_{{}_\pi} \; | \; u,v \in V\}  \; \subseteq \; V^{\otimes 2}, \quad \hbox{and}
\\
\\
&A_{{}_{\pi,\l}} \; :=\;  \frac{TV}{(R_{{}_{\pi,\l}} )}.
\end{align*}
 
\subsubsection{}   
Fix a basis $\{x,y,z\}$ for $V$. Let $f \in S^3V$. The formulas
$$
\{x,y\}_{\!{}_f} := \frac{\pd f}{\pd z}, \quad \{y,z\}_{\!{}_f} := \frac{\pd f}{\pd x}, \quad \{z,x\}_{\!{}_f} := \frac{\pd f}{\pd y},
$$ 
define a Poisson bracket $\{\cdot,\cdot\}_{\!{}_f}:=\{\cdot,\cdot\}_{\!{}_\pi}$ on $SV$ where $\pi$ is 
$$
 f_x \, \bigg(\frac{\pd }{\pd y} \wedge \frac{\pd }{\pd z} \,-\, \frac{\pd }{\pd z} \wedge \frac{\pd }{\pd y} \bigg)
 \;+\; f_y \, \bigg(\frac{\pd }{\pd z} \wedge \frac{\pd }{\pd x} \, - \, \frac{\pd }{\pd x} \wedge \frac{\pd }{\pd z} \bigg)
  \;+\; f_z \, \bigg(\frac{\pd }{\pd x} \wedge \frac{\pd }{\pd y} \, - \, \frac{\pd }{\pd y} \wedge \frac{\pd }{\pd x} \bigg)
$$
and $f_x$, $f_y$, $f_z$, are the partial derivatives of $f$. (The Poisson brackets  $\{\cdot,\cdot\}_{\!{}_f}$ are all the unimodular quadratic Poisson brackets on $SV$.)
In this case, we write $A_{\!{}_{f,\l}}$ for $A_{{}_{\pi,\l}}$. 
Thus, $A_{\!{}_{f,\l}} $ is $k\langle x,y,z\rangle$ modulo the relations
\begin{equation}
\label{defn.Af.lambda} 
\begin{cases}
xy-yx \;=\;  \l\widehat{f_z}, &  \\
yz-zy \;=\; \l \widehat{f_x}, & \\
zx-xz \;=\; \l \widehat{f_y.} &
\end{cases}
\end{equation}

\subsubsection{}
Although $\{\cdot,\cdot\}_{\!{}_f}$ and the isomorphism class of 
$A_{\!{}_{f,\l}} $ depend on the basis for $V$ our notation will not indicate that dependence.
Proposition \ref{prop.Af.basis.change} shows that $A_{\!{}_{\theta(f),\l}} \cong A_{\!{}_{f,\l\det(\theta)}}$ if $\theta \in \GL(V)$.

\subsection{Some consequences of the classification}

\subsubsection{} 
\label{ssect.Jw.not.domain}
Let $\sfw\in V^{\otimes 3}-\{0\}$.

Theorem \ref{thm.nf} gives an elementary, though not always easy to check, criterion to determine whether  a given $J(\sfw)$ is 3-Calabi-Yau. 
Often it is simpler to use the classification in Table \ref{main.table} to determine whether $J(\sfw)$ is 3-Calabi-Yau. 

\begin{thm}
\label{thm.nf} 
Let $\sfw \in V^{\otimes 3}-\{0\}$. Then $J(\sfw)$ is 3-Calabi-Yau if and only if $xy \ne 0$ for all non-zero degree-one elements 
$x,y \in J(\sfw)$. 
\end{thm}

\begin{cor}
\label{cor.nf}
Let $\sfw \in V^{\otimes 3}-\{0\}$. $J(\sfw)$ is 3-CY if and only if it is a domain. 
Up to isomorphism of algebras, there are exactly nine $J(\sfw)$ that are not  3-CY.
\end{cor}

The nine algebras referred to in Corollary \ref{cor.nf} are the $J(\sfw)$'s that are not domains.
Three of those nine algebras appear in \S\ref{ssect.Afl.not.ASR}. The other six are those labelled by a $*$   
in Table \ref{table.symm.relns.for.Jw}. The latter six are the 
non-3-CY $J(\sfw)$'s of the form $TV/(R)$ with $R \cap \Sym^2(V) \ne \{0\}$; in those cases $R \subseteq\Sym^2(V)$.

\subsubsection{} 
Every $A_{f,\l}$ is isomorphic to some  $J(\sfw)$. In fact, $A_{f,\l} =J(\sfw_0-\l \widehat{f}\,)$ provided   $\{\cdot,\cdot\}_f$ and $\sfw_0$ are defined with respect to the same basis for $V$. 

If  $J(\sfw)$ is isomorphic to some $A_{\pi,\l}$, then it is isomorphic to $A_{\overline{\sfw},-\mu(\sfw)^{-1}}$. The existence of 
the isomorphism does not depend on the choice of basis for $V$ provided the same basis is used to define $\{\cdot,\cdot\}_{{}_{\overline{\sfw}}}$ and $\sfw_0$, and hence $\mu(\sfw)$ (cf. \S\ref{ssect.cancel}).  

\subsubsection{} 
An algebra $J(\sfw)$ is isomorphic to some $A_{f,\l}$ if and only if $c(\sfw) \ne s(\sfw)$, 
and in that case $f=\overline{\sfw}$; these are the $J(\sfw)$'s in the right-most column of Table \ref{main.table}. 

\subsubsection{} 
\label{ssect.symm.relns}
Up to isomorphism of algebras, the 3-CY algebras in the  $c(\sfw) = s(\sfw)$ column of Table \ref{main.table}  form a 1-parameter family together with three additional algebras. All algebras in that column are  of the form
\begin{equation}
\label{symm.relns}
 \frac{k \langle x,y,z\rangle}{(\widehat{f_x},\,\widehat{f_y}, \, \widehat{f_z})}
\end{equation}
where  $f =\overline{\sfw}$.  Up to isomorphism of schemes, the cubic divisors on $\PP^2$ form a 1-parameter family, 
the smooth curves, together with eight singular curves. The algebra in (\ref{symm.relns}) is not 3-CY when the zero locus of 
$f$ is $\vert\!\vert\!\vert\,$, $\leftrightline \!\!\!\!\!\!\vert\!\vert\,\,$, $\varhexstar\,$, $ \vert\!  \largecircle$, or $\bigcurlywedge$.
When the zero locus of $f$ is $\largetriangleup\,$, or $\varnothing\,$, or $\propto$, respectively, the algebra in (\ref{symm.relns}) is 3-CY and isomorphic to 
\begin{itemize}
 \item 
 $J(xyz+zyx)$ with relations $yz+zy=zx+xz=xy+yx=0$, or
  \item 
   $J(xyz+zyx+\frac{1}{3}x^3)$ with relations $yz+zy+x^2=zx+xz=xy+yx=0$, or 
  \item 
   $J\big(xyz+zyx+\frac{1}{3}(x^3+y^3)\big)$ with relations $yz+zy+x^2=zx+xz+y^2=xy+yx=0$,
\end{itemize}
respectively. 

The 1-parameter family of 3-CY algebras of the form (\ref{symm.relns}) are the algebras 
$J\big(xyz+zyx+\frac{1}{3}\b(x^3+y^3+z^3)\big)$, $\b^3 \ne 1$. Such an algebra is isomorphic to  
$k\langle x,y,z\rangle$ modulo the relations
\begin{equation}
\label{1-param.family}
\begin{cases}
yz+zy+\b x^2=0, & \\ 
zx+xz+\b y^2=0, &  \\
xy+yx+\b z^2=0. & 
\end{cases}
\end{equation}
If $\b^3=1$, then $j(E)=0$ and $J(\sfw)\cong k\langle x,y,z\rangle/(x^2,y^2,z^2)$ \cite[Thm.1.1]{Sm}.  
 
Since $J(\sfw)$ only depends on $\overline{\sfw}$ when $c(\sfw)=s(\sfw)$, each algebra arises from many different $\sfw$.

\begin{thm}
\textnormal{(Theorem \ref{thm.central.squares})}
If $c(\sfw)=s(\sfw)$ and $J(\sfw)$ is a 3-CY algebra, then there is a basis $x,y,z$ for $V$ such that $J(\sfw)$ is a
 Clifford algebra of rank 8 over the central subalgebra $k[x^2,y^2,z^2]$.  
 \end{thm}

\subsubsection{}
\label{ssect.Afl.not.ASR}
If $A_{\!{}_{f,\l}}$ is not 3-Calabi-Yau it is isomorphic to  
$$
 \frac{k \langle x,y,z\rangle}{(xy, \,yz, \, zx)}
 \quad \hbox{or} \quad
  \frac{k \langle x,y,z\rangle}{(xy, \,yz-x^2, \, zx)}
 \quad \hbox{or} \quad
  \frac{k \langle x,y,z\rangle}{(xy, \,yz-x^2, \, zx-y^2)}\,
 $$
  by Corollary \ref{cor.Af.ASR}.  Since $A_{\!{}_{f,\l}}$ is  3-Calabi-Yau if and only if it is Artin-Schelter regular of dimension 3, 
  $A_{\!{}_{f,\l}}$ is  3-Calabi-Yau if and only if it is a domain.

\subsection{}
The entry in the 33-position of Table \ref{main.table} can be stated as follows: if $c(\sfw) \ne s(\sfw)$ and $E$ is $\largetriangleup$, $\varnothing\,$, or $\propto$, then 
$J(\sfw)$ fails to be 3-CY if and only if $V$ has a basis $\{x,y,z\}$ such that  $\partial_z\sfw=xy$ (or, if and only if $xy$ is zero in $J(\sfw)$).

 \subsection{Acknowledgements}
 The first author thanks the University of Washington for its hospitality during the period that this work was done. 
 Both authors are very grateful to an anonymous referee who read an earlier version of this paper very carefully, 
 spotted several mistakes, and provided numerous suggestions and feedback that have been incorporated into the final version. Thank you!

\section{Preliminaries}
\label{sect.prelim}

Always,  $k$ denotes an algebraically closed field of characteristic not 2 or 3. On rare occasions we prove 
something involving $m!$ or $S_m$ and when that happens we usually need to assume that $\fchar(k)$ does
not divide $m!$. It will be obvious to the reader when this happens so we leave the reader to fill in that additional hypothesis.
All vector spaces and algebras are defined over $k$. 
We write $\otimes$ for $\otimes_k$.

\subsection{Notation}
Let $V$ be a finite dimensional vector space.

We make the following definitions:
\begin{align*}
TV\; : = \; & \hbox{the tensor algebra} \\
\;  = \; & k \oplus V \oplus V^{\otimes 2} \oplus \cdots \\
S_m \; :=\; & \hbox{the symmetric group on $m$ letters acting on $V^{\otimes m}$ by } \\ 
&     \s(v_1 \otimes \ldots \otimes v_m)= v_{\s 1} \otimes \ldots \otimes v_{\s m}
\\
\Sym^m(V)   \;  := \; & \hbox{the symmetric tensors in $V^{\otimes m}$} \\
 \;  = \; & \{\sfw \in V^{\otimes m} \; | \;  \s(\sfw)=\sfw \hbox{ for all } \s \in S_m\}  \\
\Alt^m(V) \;  := \; & \hbox{the alternating or skew-symmetric tensors in $V^{\otimes m}$} \\
 \;  = \; & \{\sfw \in V^{\otimes m} \; | \; \s(\sfw)=\sgn(\s) \sfw \hbox{ for all } \s \in S_m\}  
\\
V^{\otimes 2} \; = \; & \Sym^2( V) \oplus \Alt^2( V) = \hbox{(symmetric)$\, \oplus\,$(skew-symmetric) tensors} \\
\phi:TV  \longrightarrow & \,TV, \quad \phi(\sfw \otimes v):=v \otimes \sfw, \; \hbox{for all $\sfw \in TV$ and $v \in V$} \\
c:V^{\otimes m}   \longrightarrow & \, V^{\otimes m}, \quad c(\sfw) :=\frac{1}{m} \sum_{j=0}^{m-1} \phi^j(\sfw), \\
s:V^{\otimes m}   \longrightarrow & \, V^{\otimes m}, \quad s(\sfw) :=\frac{1}{m!} \sum_{\s \in S_m} \s(\sfw), \\
 \end{align*}
 \begin{align*}
SV \; : = \; & \hbox{the symmetric algebra} \\
\;  = \; &  k \oplus V \oplus S^2V \oplus \cdots \\
V^{\otimes m}   \longrightarrow & \, S^mV, \; \sfw \mapsto \overline{\sfw},\; \hbox{the canonical map} \\
S^mV  \longrightarrow  & \Sym^mV,  \; f \mapsto \widehat{f}   \\
& \;  \hbox{the symmetrization map,}\\
& \;  \hbox{inverse to the map $\Sym^mV \to S^mV$, $\sfw \mapsto \overline{\sfw}$}.
\end{align*}

\subsection{The element $\sfw_0$ and the map $\mu:V^{\otimes 3} \to k$}
\label{ssect.w0}

\subsubsection{The idempotents $c$ and $s$ in the group algebra}
Let $kS_3$ denote the group algebra of $S_3$. We write ${\bf 1}$, {\bf sgn}, and ${\bf 2}$, for the trivial representation, the 
sign representation, and the 2-dimensional irreducible representation, of $kS_3$, respectively.

When $m=3$, the elements $c$ and $s$ in $kS_3$ are $c=\frac{1}{3}(1+(123)+(321))$ and 
$$
s= \hbox{$\frac{1}{2}$}c\cdot(1-(12))=  \hbox{$\frac{1}{2}$}c\cdot(1-(13))= \hbox{$\frac{1}{2}$}c\cdot(1-(23)).
$$
Since $sg=s$ for all $g \in S_3$, $s^2=s$ and $cs=sc=s$.  Since $(123)c=c$, $c^2=c$. The elements
$$
1-c, \qquad c-s, \qquad s,
$$
form a complete set of mutually orthogonal central idempotents in $kS_3$. If $M$ is a left $kS_3$-module, then
\begin{align*}
sM \; = \; & \hbox{the sum of all submodules of $M$ that are isomorphic to {\bf 1}},
\\
(c-s)M  \; = \; & \hbox{the sum of all submodules of $M$ that are isomorphic to {\bf sgn}},
\\
(1-c)M  \; = \; & \hbox{the sum of all submodules of $M$ that are isomorphic to ${\bf 2}$},
\\
cM  \; = \; & s M \oplus (c-s)M,
\\
m \; = \; & s(m)+(c-s)(m)+(1-c)(m) \hbox{ for all } m \in M.
\end{align*}
The last equality follows from the fact that $s+(c-s)+(1-c)=1$.
The second to last equality is proved by observing that $cM \subseteq s M \oplus (c-s)M$ because $c=s+(c-s)$, and 
$s M \oplus (c-s)M  \subseteq cM$ because $s=cs$ and $c-s=c(c-s)$.

\subsubsection{The element $\sfw_0$ and the $S_3$-module decomposition of $V^{\otimes 3}$}
 Let $\{x,y,z\}$ be a basis for $V$. The element
\begin{equation}
\label{defn.w0}
\sfw_0 :   \; = \; 2(c-s)(xyz) \; = \; \,\hbox{$\frac{1}{3}$}(xyz+yzx+zxy) \; - \;  \,\hbox{$\frac{1}{3}$}(zyx+xzy+yxz) \; \in \; V^{\otimes 3}
\end{equation}
is a basis for $\Alt^3(V)$. 
It follows from Lemma \ref{lem.w0} below that a change of basis for $V$ changes $\sfw_0$ by a non-zero scalar multiple.

Since $\pd_x(\sfw_0) = yz-zy$, $\pd_y(\sfw_0) = zx-xz$, and $\pd_z(\sfw_0) = xy-yx$, $J(\sfw_0) = SV$.

\begin{lem}
\label{lem.w0}
Suppose $\dim_k(V)=3$. Let $\{x,y,z\}$ be a basis for $V$, and let  $\sfw_0$ be the element in (\ref{defn.w0}). 
Then 
$$
\Alt^3(V)=k\sfw_0 =(c-s)V^{\otimes 3} \qquad \hbox{and} \qquad cV^{\otimes 3}=k\sfw_0 \oplus \Sym^3(V).
$$ 
\end{lem}
\begin{pf}
We  write $v_1v_2v_3$ for $v_1 \otimes v_2 \otimes v_3$ whenever $v_1,v_2,v_3 \in V$.

Because $cM=(c-s)M \oplus sM$ for all $kS_3$-modules $M$, $cV^{\otimes 3} =\Alt^3(V) \oplus sV^{\otimes 3} =
\Alt^3(V) \oplus \Sym^3(V)$. It remains to prove that   $\Alt^3(V)=k\sfw_0$. 

Certainly, $\sfw_0 \in \Alt^3(V)$ because  
$$
(c-s)(xyz)=\hbox{$\frac{1}{2}$}c\circ\big(1-(13)\big)(xyz)=\hbox{$\frac{1}{2}$}c(xyz -zyx)=\hbox{$\frac{1}{2}$}\sfw_0.
$$
The result of applying $c-s$ to either $yzx$, $zxy$, $zyx$, $yxz$, or $xzy$, is $\pm \sfw_0$.  
On the other hand, $(c-s)(x^2y)= \frac{1}{2}c\circ\big(1-(12)\big)(x^2y)= \frac{1}{2}c(x^2y-x^2y)=0$.
A similar calculation shows that   $(c-s)(v_1v_2v_3)=0$ whenever $\{v_1,v_2,v_3\} \subseteq \{x,y,z\}$ and $|\{v_1,v_2,v_3\} | \le 2$. To summarize, if $w$ is a word of length 3 in the letters $x,y,z$, then $(c-s)(w)$ is a multiple of $\sfw_0$. Therefore 
$(c-s)V^{\otimes 3}=k\sfw_0$. 
\end{pf}

\subsubsection{The map $\mu$}
Since $(c-s)V^{\otimes 3}=k\sfw_0$, given $\sfw \in V^{\otimes 3}$, there is a unique scalar $\mu(\sfw)$ such that 
$$
c(\sfw) \; = \; \mu(\sfw) \sfw_0 \, + \, s(\sfw).
$$
This determines a linear map $\mu:V^{\otimes 3} \to k$ that depends up to a non-zero scalar factor on the choice of basis for $V$.
It is easy to compute $\mu(\sfw)$ by using the following facts:  $\mu(\sfw)=\mu(\phi(\sfw))=\mu(c(\sfw))$, $\mu(xyz)=\frac{1}{2}$,
$\mu(zyx)= -\frac{1}{2}$, and $\mu(s(\sfw))=\mu(x^3)=\mu(xyx)=0$.

\subsubsection{}
We also note that  $V^{\otimes 3} \; = \; \Sym^3(V) \oplus \big(\Alt^2(V) \otimes V+V \otimes \Alt^2(V) \big)$ and 
$\widehat{\overline{\sfw}}=s(\sfw)$ for all $\sfw \in V^{\otimes 3}$.

 \subsection{Cyclic derivatives of elements in $TV$}

Let $x_1,\ldots,x_n$ be a basis for $V$. 

Having fixed a basis for $V$, a {\sf word} in $TV$ is an element belonging to the multiplicative semigroup generated by  $x_1,\ldots,x_n$. We 
call $x_1,\ldots,x_n$ {\sf letters}.  

Let $\ve_{x_i}$ and $\ve'_{x_i}$ be the linear maps $TV \to TV$ such that  for each word $\sfw$
\begin{equation}
\label{defn.ve.ve'}
\ve_{x_i}(\sfw ):=\begin{cases} 
					\sfv & \text{if $\sfw =x_i\sfv$,}
					\\
					0 & \text{otherwise,}
			\end{cases}	
\qquad \hbox{and} \qquad 
\ve'_{x_i}(\sfw ):=\begin{cases} 
					\sfu & \text{if $\sfw =\sfu x_i$,}
					\\
					0 & \text{otherwise.}
			\end{cases}				
\end{equation}
 
We write ${\sf span}\{a,b,\ldots\}$ for the $k$-linear span of elements $a,b,\ldots,$ in a $k$-vector space.

\begin{lem} 
\label{lem.pd.cpd}
\label{eq.2.pds}  
Let $\sfw \in V^{\otimes m}$ and $\sfu \in TV$.
\begin{enumerate}
  \item 
$\pd_{x_i}(\sfw )\; = \; m\,\ve_{x_i}(c(\sfw)) \; = \; m\,\ve'_{x_i}(c(\sfw)) $ for all $i$.
  \item 
${\sf span}\{\ve_{x_1} ( \sfu),  \ldots, \ve_{x_n} (\sfu)\}$ does not depend on the choice  of basis for $V$.
  \item 
${\sf span }\{ \pd_{x_1}( \sfw),  \ldots,   \pd_{x_n} (\sfw)\}$ does not depend on the choice  of basis for $V$.
\end{enumerate}
\end{lem}
\begin{pf}
(1)
The proofs for $\ve'_{x_i}$ are ``the same'' as those for  $\ve_{x_i}$ so we omit them.

The claim is true if $\sfw=a_1\ldots a_m$ and each $a_j \in \{x_1,\ldots,x_n\}$ because 
\begin{align*}
m\,\ve_{x_i}(c(\sfw)) \; = \; & \ve_{x_i}\big(a_1\ldots a_m \; + \; a_2\ldots a_ma_1 \; + \;  \cdots  \; + \;  a_ma_1\ldots a_{m-1}\big)
\\
\;= \; & \d_{x_i,a_1} a_2 \ldots a_m  \; + \;  \d_{x_i,a_2} a_3\ldots a_ma_1 \; + \; \cdots  \; + \;   \d_{x_i,a_m} a_1\ldots a_{m-1}
 \\
 \;= \; & \sum_{\sfw=\sfu x_i \sfv} \sfv\sfu
 \\
 \;= \; &  \pd_{x_i}(\sfw ).
 \end{align*}
Since $c$, $\ve_{x_i}$ and $\pd_{x_i}$, are linear maps, the lemma holds for all $\sfw  \in  V^{\otimes m}$.

(2)
Let $\sfu \in V^{\otimes m}$. If $\psi_i \in V^*$ is defined by $\psi_i(x_j)=\d_{ij}$, then 
$$
\ve_{x_i}(\sfu)=\big(\psi_i \otimes \id_V^{\otimes (m-1)}\big)(\sfu).
$$
Since $\psi_1,\ldots,\psi_n$ is a basis for $V^*$ it follows that
$$
{\sf span }\{\ve_{x_1} ( \sfu),  \ldots, \ve_{x_n} (\sfu)\} = {\sf span }\big\{\big(\psi \otimes \id_V^{\otimes (m-1)}\big)(\sfu) 
\; \big\vert \; \psi \in V^*\big\}.
$$
Thus, (2) is true for all $\sfu \in V^{\otimes m}$. Since $\ve_{x_i}$ is a linear map (2) holds for all $\sfu \in TV$.

(3) 
This follows from (1) and (2). 
\end{pf}

\begin{lem}
\label{lem.hat}
Let $f=v_1\ldots v_m \in S^mV$ where each $v_i \in V$. Then 
$$
\widehat{f} \; = \; \hbox{$\frac{1}{m}$}\Big(v_1 \otimes (v_{2}\ldots v_{m})\widehat{\phantom{x}} \, + \, 
 v_2 \otimes  (v_{3}  \ldots  v_{m} v_{1})\widehat{\phantom{x}}  \, + \,  \cdots
 \, + \,  v_m \otimes(v_1 \ldots  v_{m-1})\widehat{\phantom{x}}  \Big).
 $$
 \end{lem}
 \begin{pf}
For $i=1, \dots, m$, let $G_i:=\{\s\in S_m\; | \; \s(i)=i\}$ and $G_i' := \{\s \in S_m \; | \; \s(1)=i\}$.  Since $S_m$ is the disjoint union of the sets $G_i'$, 
%  = \{\s \in S_m \; | \; \s(1)=i\}$, 
$1 \le i \le m$,
\begin{align*}
\widehat{f} \; = \; & \frac{1}{m!} \sum_{\s \in S_m} v_{\s 1}\otimes \ldots \otimes v_{\s m}
\\
\; = \; & \frac{1}{m!}\Bigg(\sum_{\s \in G_1'} v_{\s 1}\otimes \ldots \otimes v_{\s m} +
 \sum_{\s \in G_2'} v_{\s 1}\otimes \ldots \otimes v_{\s m} + \cdots 
 \\
 & \qquad \qquad \qquad \cdots +  \sum_{\s \in G_m'} v_{\s 1}\otimes \ldots \otimes v_{\s 
m} \Bigg)
\\
 \; = \; & \frac{1}{m!}\Bigg(v_1 \otimes \sum_{\s \in G_1} v_{\s 2}\otimes \ldots \otimes v_{\s m} +
 v_2 \otimes \sum_{\s \in G_2} v_{\s 3}\otimes \ldots \otimes v_{\s m}  \otimes v_{\s 1} + \cdots 
 \\
 & \qquad \qquad \qquad \cdots +  v_m \otimes \sum_{\s \in G_m} v_{\s 1}\otimes \ldots \otimes v_{\s (m-1)} \Bigg)
\\
\; = \; &   \frac{1}{m}\Big(v_1 \otimes (v_{2}\ldots v_{m})\widehat{\phantom{x}}+
 v_2 \otimes  (v_{3}  \ldots  v_{m} v_{1})\widehat{\phantom{x}} + \cdots
 +  v_m \otimes(v_1 \ldots  v_{m-1})\widehat{\phantom{x}}  \Big),
 \end{align*}
 as claimed.
 \end{pf}
 
The next result is a special case of Lemma \ref{lem.hat}. 

\begin{lem}
\label{lem.g.hat}
Let $\{x_1,\ldots,x_n\}$ be a basis for $V$. Let $g \in S^2V$. Let $g_i$ denote the partial derivative of $g$ with respect
to $x_i$. 
If we view $g_i$ as an element in $TV$, then 
$$
\widehat{g}\; = \; \frac{1}{2} \sum_{i=1}^n x_i g_i  \; = \; \frac{1}{2} \sum_{i=1}^n  g_i x_i.
$$
\end{lem}

\begin{lem}
\label{lem.pd.hat}
Let $\{x_1,\ldots,x_n\}$ be a basis for $V$. If $f \in S^mV$, then
\begin{equation}
\label{lem.hat.circ.d}
 \partial_{x_i}(\widehat{f}\,) \; = \; m \, \ve_{x_i}\big(\widehat{f}\,\big) \; = \; m \, \ve'_{x_i}\big(\widehat{f}\,\big)   \; = \;   \widehat{\Bigg(\frac{\pd f}{\pd {x_i}}\Bigg)} 
\end{equation}
for all  $i$.
\end{lem}
\begin{pf}
Because $\pd_{x_i}$, $\ve_{x_i}$, $\frac{\pd}{\pd x_i}$, and $f \mapsto \widehat{f}$, are linear maps it suffices to prove the 
result when $f=v_1\ldots v_m$ where each $v_j \in \{x_1,\ldots,x_n\}$. We therefore assume that $f$ is of this form.

Let $x \in \{x_1,\ldots,x_n\}$. 
It follows from the expression for $\widehat{f}$ in Lemma \ref{lem.hat} that 
\begin{align*}
m\ve_x(\widehat{f})\; = \; &  \ve_x\Big(v_1 \otimes (v_{2}\ldots v_{m})\widehat{\phantom{x}}+
 v_2 \otimes  (v_{3}  \ldots  v_{m} v_{1})\widehat{\phantom{x}} + \cdots
 +  v_m \otimes(v_1 \ldots  v_{m-1})\widehat{\phantom{x}}  \Big)
 \\
 \; = \; &  \d_{xv_1} (v_{2}\ldots v_{m})\widehat{\phantom{x}}+
\d_{xv_2} (v_{3}  \ldots  v_{m} v_{1})\widehat{\phantom{x}} + \cdots
 + \d_{xv_m} (v_1 \ldots  v_{m-1})\widehat{\phantom{x}} 
  \\
 \; = \; & ( \d_{xv_1} v_{2}\ldots v_{m}+
\d_{xv_2} v_{3}  \ldots  v_{m} v_{1}+ \cdots
 + \d_{xv_m}v_1 \ldots  v_{m-1} )\widehat{\phantom{x}}
   \\
 \; = \; & \widehat{\Bigg(\frac{\pd f}{\pd x}\Bigg)}.
 \end{align*}
By Lemma \ref{eq.2.pds}(1), $\partial_{x_i}(\hat f)  = m \, \ve_{x_i}(c(\widehat{f}))=m\ve_{x_i}(\widehat{f})$.
 The proof is complete. 
\end{pf}
 
By Lemma \ref{lem.pd.cpd}(3), $J(\sfw)$ depends only on $\sfw$ and not on the choice of basis for $V$.

\subsection{} 
We are now ready to prove that $J(\sfw)$ is 3-Calabi-Yau if and only if it is Artin-Schelter regular of dimension 3. 
 
\begin{prop}
\label{prop.Jw.xMx}
Let $\{x_1,\ldots,x_n\}$ be a basis for $V$. Let $\sfx^\sT=(x_1,\ldots,x_n)$. 
Let  $\sfw \in V^{\otimes m}$. If $\sM$ is the unique matrix  
such that $c(\sfw)=\sfx^\sT \sM\sfx$, then
\begin{equation}
\label{eq.Jw.xMx}
\sM\sfx \; = \; \frac{1}{m} \, \begin{pmatrix}   \pd_{x_1}\!( \sfw) \\ \\ \vdots \\ \\   \pd_{x_n}\!( \sfw)  \end{pmatrix} 
\;= \; (\sfx^\sT  \sM)^\sT.
\end{equation}
\end{prop}
\begin{pf}
If $m \ge 1$ and $\sfw \in V^{\otimes m}$, then 
$$
\sfw \; = \; \sum_{i=1}^n x_i\ve_{x_i}\!(\sfw) 
\;=\;  \begin{pmatrix}x_1 & \cdots & x_n  \end{pmatrix} 
 \begin{pmatrix} \ve_{x_1}\!( \sfw) \\ \\ \vdots \\ \\   \ve_{x_n}\!(\sfw)  \end{pmatrix} .
 $$
When $\sfw=c(\sfw)$, this equality in conjunction with Lemma \ref{lem.pd.cpd}  gives
 $$
\sfx^\sT \sM\sfx \; = \; c(\sfw)  \; = \;\sfx^\sT 
 \begin{pmatrix} \ve_{x_1}\!(c(\sfw)) \\ \\ \vdots \\ \\   \ve_{x_n}\!(c(\sfw))  \end{pmatrix}  \; = \; \hbox{$\frac{1}{m}$}\,\sfx^\sT 
  \begin{pmatrix}   \pd_{x_1}( \sfw) \\ \\ \vdots \\ \\   \pd_{x_n}( \sfw)  \end{pmatrix}.
 $$
 Because the entries in $\sfx^\sT$ are a basis for $V$ we can cancel the $\sfx^\sT$ factors in the previous equality. 
 Doing that gives the first equality in (\ref{eq.Jw.xMx}). 
The second equality is proved in a similar way by using $\ve_{x_i}'$ in place of $\ve_{x_i}$. 
\end{pf}

\begin{cor} 
\label{cor.Jw.standard}
Suppose $\dim_k(V)=3$ and let $\sfw \in V^{\otimes 3}-\{0\}$. If the subspace $R_\sfw\subseteq V^{\otimes 2}$ spanned by
the cyclic partial derivatives has dimension 3, then $J(\sfw)$ is standard in the sense of \S\ref{ssect.standard} and the matrix $Q$ in 
\S\ref{ssect.standard} is the identity.
\end{cor}
\begin{pf}  
Fix a basis $\{x_1,x_2,x_3\}$ for $V$ and adopt the notation in Proposition \ref{prop.Jw.xMx}. 
By definition, $J(\sfw)$ is $k\langle x_1,x_2,x_3\rangle$ modulo the ideal generated by the entries in $\sM\sfx$.
By hypothesis, the entries in $\sM\sfx$ are linearly independent. Since $\sM\sfx = (\sfx^\sT  \sM)^\sT$, $J(\sfw)$ is standard.
Furthermore, the matrix $Q$ with the property  $\sfx^\sT\sM=(Q\sM\sfx)^\sT$ is the identity.   
\end{pf}

\begin{prop}
\label{prop.3ASR=3CY}
Let $V$ be a 3-dimensional vector space, and $R$ a 3-dimensional subspace of $V^{\otimes 2}$. The algebra $TV/(R)$ 
is a twisted 3-Calabi-Yau algebra if and only if it is a 3-dimensional Artin-Schelter regular algebra.\footnote{The ``same proof'' shows that if
$V$ and $R \subseteq V^{\otimes 3}$ have dimension 2, then $TV/(R)$ 
is a twisted 3-Calabi-Yau algebra if and only if it is a 3-dimensional cubic Artin-Schelter regular algebra.}
\end{prop}
\begin{pf}
By  \cite[Lem. 1.2]{RRZ}, a connected graded algebra $A$ is twisted $d-$Calabi-Yau if and only if $\gldim(A)=d<\infty$ and 
(\ref{AS-Gor-cond}) holds for some $\ell \in \ZZ$. Thus, 3-dimensional Artin-Schelter regular algebras are twisted 3-Calabi-Yau.

Let $A=TV/(R)$.
If $V$ is any graded vector space and $R$ any graded subspace of $TV$ that is zero in degrees 0 and 1, the minimal projective resolution of the trivial module ${}_Ak$  begins $A \otimes R \to A \otimes V \to A \to k \to 0$. If such an $A$ is twisted 3-Calabi-Yau, then it has global dimension 3 so the full minimal resolution of ${}_Ak$ is 
$0 \to A^n \to A \otimes R \to A \otimes V \to A \to k \to 0$ for some $n$. 
Because this resolution is minimal, the dimension of $\Ext_A^3(k,A)$ is $n$. 
Hence by \cite[Lem. 4.3(a)]{YZ},  the minimal resolution is 
\begin{equation}
\label{eq.minl.res}
0 \to A(-\ell) \to A \otimes R \to A \otimes V \to A \to k \to 0.
\end{equation}
However, (\ref{AS-Gor-cond}) implies that $\Hom_A(-,A)$ applied to (\ref{eq.minl.res}) gives a deleted resolution of $k(\ell)$. Because $R$ is concentrated in degree 2, $\ell=3$; i.e., the minimal resolution is $0 \to A(-3)  \to A(-2)^3  \to A(-1)^3  \to A   \to k \to 0$. Hence the Hilbert series of $A$ is $(1-t)^{-3}$. Thus, $\GKdim(A)=3$ and 
$A$ is Artin-Schelter regular of dimension 3.  
 \end{pf}

 \begin{cor}
\label{cor.3ASR=3CY}
Let $V$ be a 3-dimensional vector space and $\sfw \in V^{\otimes 2}-\{0\}$. Then $J(\sfw)$ is 3-Calabi-Yau if and only if it is  
 Artin-Schelter regular of dimension 3. Thus, if  $J(\sfw)$ is 3-Calabi-Yau, it is a domain.
\end{cor}
\begin{pf}
Suppose $J(\sfw)$ is Artin-Schelter regular of dimension 3. 
 By \cite[Cor. 9.3]{VdB3}, the Nakayama automorphism of a  3-dimensional quadratic Artin-Schelter regular algebra  
 is induced by $\sfx \mapsto Q^{-\sT}\sfx$. But $Q=I$ by Corollary \ref{cor.Jw.standard}, so $J(\sfw)$ is 3-Calabi-Yau. Proposition \ref{prop.3ASR=3CY} proves the converse. By \cite[Thm. 3.9]{ATV2}, 3-dimensional Artin-Schelter regular algebras are domains. 
 \end{pf}

\subsection{The algebras $A_{{}_{\pi,\l}}$ and $A_{{}_{f,\l}}$ }

\begin{lem}
\label{lem.Alt2.A.lambda}
Suppose $\dim_k(V)=3$. Let $\l \in k$,
and let $\{\cdot,\cdot\}_{{}_\pi}$ be a quadratic Poisson bracket  on $SV$ with associated bivector field $\pi$.  
If $R\subseteq V^{\otimes 2}$ is a subspace such that $TV/(R) \cong A_{{}_{\pi,\l}}$ as graded $k$-algebras, then 
$R\cap \Sym^2(V)=\{0\}$.  
\end{lem}
\begin{pf}
Suppose $\Theta:TV/(R) \to A_{{}_{\pi,\l}}$ is a graded $k$-algebra isomorphism.
Let $\theta$ be the restriction of $\Theta$ to $V$. 
Then $(\theta\otimes \theta)(R)=R_{{}_{\pi,\l}}$.  Since $(\theta \otimes \theta)(\Sym^2V) = \Sym^2V$,  
$$
(\theta\otimes \theta)(R\cap \Sym^2V)=R_{{}_{\pi,\l}} \cap \Sym^2V.
$$
Equivalently,  $R\cap \Sym^2V=(\theta\otimes \theta)^{-1} \big( R_{{}_{\pi,\l}} \cap \Sym^2V  \big)$.

Let $r \in R_{{}_{\pi,\l}}  \cap \Sym^2V$. 
Since $r \in R_{{}_{\pi,\l}}  $,  $r= u \otimes v - v \otimes u -\l \widehat{ \{u,v\}}$ for some $u,v \in V$. Since $r$ and $ \widehat{ \{u,v\}}$
are in $\Sym^2V$, so is $u \otimes v - v \otimes u$. But $\Alt^2V \cap \Sym^2V=\{0\}$ so $u \otimes v - v \otimes u=0$; 
hence $u$ is a multiple of $v$, or vice versa;
but $\{u,u\}=\{v,v\}=0$ so $\{u,v\}=0$; hence $r=0$. Thus, $ R_{{}_{\pi,\l}}  \cap \, \Sym^2V=\{0\}$.
Hence $R\cap \Sym^2(V)=\{0\}$. 
\end{pf}

\begin{prop}
\label{prop.Af.basis.change}
Let $f \in S^3V$ and $\l \in k$. Every $\theta \in \GL(V)$ extends to an isomorphism 
$A_{\!{}_{f,\l\det(\theta)}} \longrightarrow A_{\!{}_{\theta(f),\l}}$  of graded $k$-algebras.
\end{prop}
\begin{pf}
The notation $\theta(f)$ in the statement of the proposition  implicitly refers to the unique algebra automorphism $SV \to SV$ whose restriction to $V$ is $\theta$; i.e., $\theta$ also denotes that unique extension. Likewise, we use the letter $\theta$ to denote
the unique  automorphism $TV \to TV$ whose restriction to $V$ is $\theta$.  

We write $R_{\!{}_{f,\l\det(\theta)}}$ for the kernel of the map $V^{\otimes 2} \to A_{\!{}_{f,\l\det(\theta)}}$, and
$R_{\!{}_{\theta(f),\l}}$ for the kernel of the map $V^{\otimes 2} \to A_{\!{}_{\theta(f),\l}}$.
To prove the proposition we must show that  
\begin{equation}
\label{eq.2.relns}
\theta\big(R_{\!{}_{f,\l\det(\theta)}}\big) \; = \; R_{\!{}_{\theta(f),\l}}\,.
\end{equation}
Fix a basis $\{x,y,z\}$ for $V$. Elements in $R_{\!{}_{f,\l\det(\theta)}}$ are of the form
$$
[u,v]\, - \, \l \det(\theta)\Big( \{u,v\}_{\!{}_f}\Big)\widehat{\phantom{\Big)}}
 $$
where $u,v \in V$. 
The calculation
\begin{align*}
\theta \Bigg([u,v] - \l \det(\theta)\Big( \{u,v\}_{\!{}_f}\Big)\widehat{\phantom{\Big)}} \, \Bigg) 
& \;=\;
\theta \Bigg([u,v] - \l \det(\theta)\bigg( \frac{\sfd u \wedge \sfd v \wedge \sfd f}{\sfd x \wedge \sfd y \wedge \sfd z}\bigg)\widehat{\phantom{\bigg)}}\, \Bigg) 
\\
& \;=\;
[\theta(u),\theta(v)] -   \l \det(\theta)\bigg(\frac{\sfd \theta(u) \wedge \sfd \theta(v) \wedge \sfd \theta(f)}{\sfd \theta(x) \wedge \sfd \theta(y) \wedge \sfd \theta(z)}\,\bigg)\widehat{\phantom{\bigg)}}\,
\\
& \;=\;
[\theta(u),\theta(v)] -   \l \bigg(\frac{\sfd \theta(u) \wedge \sfd \theta(v) \wedge \sfd \theta(f)}{\sfd x \wedge \sfd y \wedge \sfd z}\,\bigg)\widehat{\phantom{\bigg)}}\,
\\
& \;=\;
[\theta(u),\theta(v)] -   \l \bigg( \Big\{ \theta(u), \theta(v)\Big\}_{\theta(f)}\,\bigg)\widehat{\phantom{\bigg)}}
\end{align*}
shows that $\theta\big(R_{\!{}_{f,\l\det(\theta)}}\big)$ is contained in $R_{\!{}_{\theta(f),\l}}$. Since $\theta$ is an automorphism of $V^{\otimes 2}$ and $R_{\!{}_{f,\l\det(\theta)}}$ and $R_{\!{}_{\theta(f),\l}}$ have the same dimension, the equality in (\ref{eq.2.relns}) holds. 
\end{pf}

\begin{lem}
\label{lem.isom.DML}
Let $\pi_1$ and $\pi_2$ be bivector fields on $V$, and $\l_1,\l_2 \in k$. 
Then $R_{{}_{\pi_1,\l_1}} =R_{{}_{\pi_2,\l_2}}$ if and only if $\l_1\pi_1=\l_2\pi_2$. 
\end{lem}
\begin{pf}
We will write $\widehat{\{\cdot,\cdot\}}_{i}$ for $\widehat{\{\cdot,\cdot\}}_{\pi_i}$ 
Clearly, $R_{{}_{\pi_1,\l_1}} =R_{{}_{\pi_2,\l_2}}$ if and only if $\l_1\widehat{\{u,v\}}_{{}_1} = \l_2\widehat{\{u,v\}}_{{}_2}$
for all $u,v \in V$. Since the function $g \mapsto \widehat{g}\,$ is an isomorphism from $S^2V$ to $\Sym^2(V)$, 
$R_{{}_{\pi_1,\l_1}} =R_{{}_{\pi_2,\l_2}}$ if and only if $\l_1\{u,v\}_{{}_1} = \l_2\{u,v\}_{{}_2}$ for all $u,v \in V$;
i.e.,  if and only if $\l_1\pi_1=\l_2\pi_2$. 
\end{pf}

\subsection{Plane cubic curves}
\label{ssect.cubics}

Our classification of 3-CY algebras uses the classification of plane cubic curves over $k$.

The scheme-theoretic zero loci of elements $f$ and $f'$ in $S^3V$ are isomorphic
if and only if $\theta(f)=f'$ for some $\theta \in \GL(V)$. ``Standard'' forms of $f$ for each isomorphism class of cubic divisors 
can be found in  \cite[Ch.4]{Hulek}, for example. We modify some of the standard forms and use the forms in Table
\ref{table.symm.relns.for.Jw}.  See  \cite[Ch.4]{Hulek} for more details.

$$
\begin{array}{|c|c|l|}
\hline
\hbox{Curve} & f & \hbox{Relations for } k\langle x,y,z\rangle/(\widehat{f_x},\,\widehat{f_y}, \, \widehat{f_z}) \phantom{\Big)}
 \\
\hline
 \vert\!\vert\!\vert   &x^3 & x^2 \; = \; 0 \phantom{\Big(xxxxxxxxxxxxxxxxxxxxxxx|xxx}  *\\
\hline
\leftrightline  \!\!\!\!\!\vert\!\vert  &x^2y &x^2 \; = \; xy+yx  \; = \; 0 \phantom{\Big(xxxxxxxxxxxxxxxxxx}  *\\
\hline
\varhexstar   &x^3+y^3 &x^2 \; = \; y^2   \; = \; 0 \phantom{\Big(xxxxxxxxxxxxxxxxxxxxxx}  * \\
\hline
\vert\!  \largecircle &x^2z+xy^2 & x^2 \; = \; xz+zx+y^2 \; = \; xy+yx  \; = \; 0 \phantom{\Big(xx|xxxx}  *\\
\hline
\bigcurlywedge &x^2z+\frac{1}{3}y^3 & x^2\; = \; y^2 \; = \;  xz+zx \; = \; 0 \;  \phantom{\Big(xx|xxxxxxxxxxx}  *\\
\hline
\largetriangleup &xyz & yz+zy \; = \; zx+xz \; = \; xy+yx \; = \; 0 \phantom{\Big(} \\
\hline
\varnothing &xyz+\frac{1}{3}x^3 &  zx+xz \; = \;  xy+yx \; = \; yz+zy+x^2  \; = \; 0  \phantom{\Big(} \\
\hline
 \propto  &xyz+\frac{1}{3}x^3+\frac{1}{3}y^3 & xy+yx \; = \; zy+yz+x^2 \; = \;  xz+zx+y^2  \; = \; 0  \phantom{\Big(} \\
\hline 
j(E)=0  &x^3+ y^3+z^3  & x^2 \; = \; y^2 \; = \;  z^2  \; = \; 0  \phantom{\Big(xxxxx|xxxxxxxxxxxx}  *\\
\hline 
j(E)\ne 0  & 2xyz + \frac{1}{3}\b( x^3+y^3+z^3)   & \b x^2 +yz+zy \; = \; \b y^2+xz+zx   \phantom{\Big(} \\
&   \b^3 \in k-\{0,1,-8\} & \phantom{\b x^2 +yz+zy} \; = \; \b z^2+xy+yx  \; = \; 0   \phantom{\Big(} \\
\hline 
\end{array}
$$
\vskip -.3in

\begin{table}[htdp]
\caption{The algebras $k\langle x,y,z\rangle/(\widehat{f_x},\,\widehat{f_y}, \, \widehat{f_z})$ up to isomorphism. }
\label{table.symm.relns.for.Jw}
\end{table}

If $c(\sfw)=s(\sfw)$ and $f=\overline{\sfw} \in S^3V$, then $J(\sfw) \cong
k \langle x,y,z\rangle/(\widehat{f_x},\,\widehat{f_y}, \, \widehat{f_z})$ by Proposition \ref{prop.Jw.relns}. The defining relations for these algebras are in the right-most column of Table \ref{table.symm.relns.for.Jw}. If $\sfw \in V^{\otimes 3}-\{0\}$, then $R_\sfw \subseteq\Sym^2(V)$ if and only if $c(\sfw)=s(\sfw)$ so the algebras 
in Table \ref{table.symm.relns.for.Jw} are the only $J(\sfw)$s for which $R_\sfw \subseteq\Sym^2(V)$.

\subsubsection{}
In the last line of Table \ref{table.symm.relns.for.Jw}, if $\b^3 \in \{0,-8\}$, then $E$ is $\largetriangleup$; if $\b^3=1$, then $j(E)=0$.

\subsubsection{}
The algebras labelled by a $*$ in Table \ref{table.symm.relns.for.Jw} are not domains and therefore  not 3-Calabi-Yau
by Corollary \ref{cor.3ASR=3CY}. Proposition \ref{prop.not.CY} shows that the other algebras in 
Table \ref{table.symm.relns.for.Jw} are 3-Calabi-Yau.

 \section{The relation between $J(\sfw)$ and $A_{\!{}_{f,\l}}$}

 \subsection{Summary of the results in this section}
 Proposition \ref{prop.Jw.relns} and Theorem \ref{thm.3.2 } show that the five statements in the left-hand column of the following table are equivalent, and 
the five statements in the right-hand column  are equivalent. 

$$ 
 \begin{array}{|l|l|} 
\hline
c(\sfw)=s(\sfw) & c(\sfw)\ne s(\sfw) \phantom{\Big)} 
\\
\hline
J(\sfw) =J(\widehat{f}\,)  \hbox{ where } f=\overline{\sfw} &  J(\sfw)\neq J(\widehat {f}) \hbox{ for any } f  \phantom{\bigg)}   \\
\hline
J(\sfw) \not\cong A_{\pi,\l}  \hbox{ for any } (\pi,\l) &   J(\sfw) \cong A_{{}_{\overline{\sfw},-\mu(\sfw)^{-1}}} \phantom{\Big)^3} 
\\
\hline
 R_\sfw \subseteq\Sym^2V   &  R_\sfw \not\subseteq \Sym^2V  \phantom{\Big)^3}  \\ 
\hline
 R_\sfw \cap \Sym^2V\ne \{0\}  &  R_\sfw \cap \Sym^2V= \{0\} \phantom{\bigg)}   \\
 \hline
\end{array}
$$

  \subsection{}
Let $\{x,y,z\}$ be a basis for $V$. By definition, $J(\sfw)$ is $k\langle x,y,z\rangle$ modulo the relations
$\pd_{x}(\sfw)=  \pd_{y}(\sfw)  =  \pd_{z}(\sfw)  =0$.

\begin{prop}
\label{prop.Jw.relns}
Let $\sfw \in V^{\otimes 3}$.
Fix a  basis $\{x,y,z\}$ for $V$. Define $\mu:V^{\otimes 3} \to k$ as in \S\ref{ssect.w0}  with respect to
$\{x,y,z\}$. Then $J(\sfw)$ is defined by the relations
\begin{align*}
\mu(\sfw)(yz-zy) \;+\; &  \widehat{\{y,z\}}_{\overline{\sfw}} \;= \; 0,
 \\
\mu(\sfw)(zx-xz) \;+\; & \widehat{\{z,x\}}_{\overline{\sfw}} \;= \; 0,
 \\
\mu(\sfw)(xy-yx) \;+\;  &   \widehat{\{x,y\}}_{\overline{\sfw}}  \;= \; 0.
\end{align*}
In particular, if $c(\sfw)=s(\sfw)$, then the defining relations for $J(\sfw)$ are 
$$ 
 \widehat{\{y,z\}}_{\overline{\sfw}} \;= \;  \widehat{\{z,x\}}_{\overline{\sfw}} \;= \;    \widehat{\{x,y\}}_{\overline{\sfw}}  \;= \; 0.
$$
\end{prop}
\begin{pf}
Let $\sfw_0$ be the element in \S\ref{ssect.w0} defined with respect to the ordered basis $\{x,y,z\}$. Thus, $c(\sfw)=\mu(\sfw)\sfw_0+s(\sfw)$.
We have  
 $$ 
   \pd_{x}(\sfw)    \;= \;   \pd_{x}(c(\sfw))  \;= \;    \pd_{x}\big(\mu(\sfw)\sfw_0+s(\sfw)\big)   \;= \;  \mu(\sfw)(yz-zy)+ 
 \Bigg(\frac{\pd \overline{s(\sfw)}}{\pd x}\Bigg)\!\widehat{\phantom{\Bigg)}} 
 $$
 where the last equality is given by Lemma \ref{lem.pd.hat}. Since $\overline{s(\sfw)}=\overline{\sfw}$,  
 \begin{equation}
 \label{relns.Jw}
 \pd_{x}(\sfw) \;=\;   \mu(\sfw)(yz-zy)+ \widehat{\frac{\pd \overline{\sfw}}{\pd x}} \; = \;   \mu(\sfw)(yz-zy) + \widehat{\{y,z\}}_{\overline{\sfw}}.
 \end{equation} 
 Similar formulae hold for  $\pd_{y}(\sfw)$ and $  \pd_{z}(\sfw)$. 
\end{pf}

\begin{thm} 
\label{thm.3.2 }
 Fix a basis $\{x,y,z\}$ for $V$ and let $\sfw_0$ be the element in \S\ref{ssect.w0}. 
Let $f \in S^3V$, $\l \in k$, and $\sfw \in V^{\otimes 3}$.  
\begin{enumerate}
  \item 
$A_{\!{}_{f,\l}} = J(\sfw_0- \l \widehat{f}\,)$.
  \item 
If $c(\sfw) \ne s(\sfw) $, then   $J(\sfw) = A_{{}_{\overline{\sfw},\nu}}$ where $\nu=-\mu(\sfw)^{-1}$.   
\item{}
If $c(\sfw) = s(\sfw)$, then $J(\sfw)$ is not isomorphic to any $A_{\!{}_{\pi,\l}}$. 
\end{enumerate}
\end{thm}
\begin{pf}
(1)
By Lemma \ref{lem.pd.hat}, $\pd_{z}(\sfw_0- \l \widehat{f})= x \otimes y - y \otimes x  -\l \widehat{f_z}$. 
There are similar expressions for $ \pd_{y} $ and $ \pd_{x}$. The result follows.
 
(2)
This is a restatement  of Proposition \ref{prop.Jw.relns} when $c(\sfw) \ne s(\sfw) $.

(3)
Suppose $c(\sfw) = s(\sfw)$.
By Proposition \ref{prop.Jw.relns}, $R_\sfw$ is spanned by the elements
$$
 \widehat{\{x,y\}}_{\overline{\sfw}}, \quad  \widehat{\{z,x\}}_{\overline{\sfw}}, \quad 
 \widehat{\{y,z\}}_{\overline{\sfw}},
 $$
which belong to $\Sym^2(V)$. By Lemma \ref{lem.Alt2.A.lambda}, the relations for $A_{\!{}_{\pi,\l}}$ are not contained in $\Sym^2(V)$.
Hence $J(\sfw)$ is not isomorphic to $A_{\!{}_{\pi,\l}}$ for any $(\pi,\l)$.
  \end{pf}
  
  \begin{cor}
  \label{cor.thm.3.2}
 For all $f \in S^3V$ and $\l \in k$, $A_{\!{}_{f,\l}}$ is standard in the sense of \S\ref{ssect.standard}.
 \end{cor}
\begin{pf}
The dimension of the kernel of the map $V^{\otimes 2} \to A_{\!{}_{f,\l}}$ is obviously 3. 
Because $A_{\!{}_{f,\l}}$ is isomorphic to $J(\sfw_0- \l \widehat{f}\,)$ it is standard by 
Corollary \ref{cor.Jw.standard}.
\end{pf}

\subsubsection{Theorem \ref{thm.3.2 } and its proof}
Although $\sfw_0$ and $\{\cdot,\cdot\}_{{}_f}$ depend on the choice of basis for $V$,
in the statement of Theorem \ref{thm.3.2 } these dependencies ``cancel out''. Whichever basis is used
for $V$,  $A_{\!{}_{f,\l}} = J(\sfw_0- \l \widehat{f}\,)$. Likewise, when $c(\sfw) \ne s(\sfw)$,   $J(\sfw) = A_{{}_{\overline{\sfw},-\mu(\sfw)^{-1}}}$ 
provided the same basis is used to define $\sfw_0$ and $\{\cdot,\cdot\}_{{}_{\overline{\sfw}}}\,$.

 \section{When is  $A_{\!{}_{f,\l}}$ 3-Calabi-Yau?}
 
 \subsection{}
 By Lemma \ref{lem.g.hat},
$A_{\!{}_{f,\l}}$ is $k\langle x,y,z\rangle$ modulo the relations  
\begin{equation}
\label{relns-Afl}
\begin{cases}
yz-zy-  \hbox{$\frac{1}{2}$}  \l \, (f_{xx}x+f_{xy}y+f_{xz}z)  \; = \; 0,  & \\
zx-xz-\hbox{$\frac{1}{2}$}   \l \, (f_{yx}x+f_{yy}y+f_{yz}z)  \; = \; 0, & \\
 xy-yx-\hbox{$\frac{1}{2}$}  \l \, (f_{zx}x+f_{zy}y+f_{zz}z)  \; = \; 0. &
 \end{cases}
\end{equation}
These relations are the entries in $\sM \sfx  \; = \; \sM(x,y,z)^\sT$  where
\begin{equation}
\label{M}
\sM\; = \; - \begin{pmatrix} 
\frac{1}{2}\l f_{xx} & \frac{1}{2}\l f_{xy}+z &\frac{1}{2}\l f_{xz}-y \phantom{\Big(} \\
\frac{1}{2}\l f_{yx}-z &\frac{1}{2}\l f_{yy} & \frac{1}{2}\l f_{yz}+x \phantom{\Big(} \\
\frac{1}{2}\l f_{zx}+y & \frac{1}{2}\l f_{zy}-x & \frac{1}{2}\l f_{zz} \phantom{\Big(}
\end{pmatrix}.
\end{equation} 

\begin{lem}
 \label{lem.22.entry-1}
Let $\{x,y,z\}$ be a basis for $V$, $\theta \in \GL(V)$, and $g \in \big\{0,\frac{1}{3}x^3, \frac{1}{3}(x^3+y^3)\big\}$.
\begin{enumerate}
  \item 
 $A_{{}_{xyz+g,\l}}$ is 3-Calabi-Yau if and only if $ \l^2 \ne 4$. 
  \item 
 $A_{{}_{\theta(xyz+g),\l}}$ is 3-Calabi-Yau if and only if $ \l^2\det(\theta)^2 \ne 4$.
\end{enumerate}
\end{lem}

\begin{pf}
(1)
Let $f=xyz+g$,  $\a=\frac{1}{2}(2-\l)$, and $\b=\frac{1}{2}(2+\l)$. Since $g_{xy}=g_{yx}=g_z=0$, it follows from (\ref{relns-Afl}) that $A_{\!{}_{f,\l}}$ is defined by the relations 

\begin{equation}
\label{22.entry.relns}
 \begin{cases}
\phantom{i} 
\a yz \, - \, \b zy\, - \, \hbox{$\frac{1}{2}$} \l \, g_{xx}x \; = \; 0, & \\
\phantom{i}  
\a zx \, - \, \b xz\, - \, \hbox{$\frac{1}{2}$} \l \, g_{yy}y   \; = \; 0, &\\
\phantom{i}  
 \a xy\, - \, \b yx  \; = \; 0. &
\end{cases}
\end{equation}

($\Rightarrow$) 
Suppose $A_{\!{}_{f,\l}}$ is 3-Calabi-Yau. Then it is a domain by Corollary \ref{cor.3ASR=3CY}.
If $\l^2$ were equal to 4, then exactly one of $\a$ and $\b$ would be 0  so $A_{\!{}_{f,\l}}$  would not be a domain; this is not the case so we conclude that $\l^2 \ne 4$.

($\Leftarrow$) 
Suppose  $\l^2 \ne 4$. Thus $\a\b\ne 0$. Let $\overline{\sM}$ denote the image of the matrix $\sM$ in (\ref{M}) in the ring of $3\times 3$ matrices over $SV$. 

By Corollary \ref{cor.thm.3.2}, $A_{\!{}_{f,\l}}$ is standard in the sense of \S\ref{ssect.standard} so, 
by Theorem \ref{thm.ATV1.1}, $A_{\!{}_{f,\l}}$ is a 3-dimensional Artin-Schelter regular algebra if and only if   
$\rank(\overline{\sM}) \ge 2$ at all points of $\PP^2$. Thus, to complete the proof of (1) it suffices to  
show that the common zero locus of the $2 \times 2$ minors of 
$$
\overline{\sM} \; =  \; - \begin{pmatrix} 
\frac{1}{2}\l g_{xx} & \b z &  -\a y \phantom{\Big(} \\
 -\a z &\frac{1}{2}\l g_{yy} &  \b x  \phantom{\Big(} \\
 \b y &  -\a x & 0 \phantom{\Big(} 
\end{pmatrix}
$$
is empty. Suppose $p=(a,b,c)\in \PP^2$ belongs to the common zero locus of the $2 \times 2$ minors of $\overline{\sM}$.
Therefore $\a \b x^2$,  $\a \b y^2$,  and $\frac{1}{4}  \l^2 \, g_{xx}g_{yy} +\a \b z^2$, vanish at $p$. In particular, $a=b=0$. 
It follows that  $g_{xx}(p)=g_{yy}(p)=0$. Hence $\a\b z^2$ also vanishes at $p$, whence $c=0$ also. Hence no such $p$ exists. 

(2)
This follows from (1) because, by Proposition \ref{prop.Af.basis.change}, $A_{\!{}_{\theta(f),\l}} \cong A_{\!{}_{f,\l\det(\theta)}}$.  
\end{pf}
 
The next result tells us when $J(\sfw)$ is 3-Calabi-Yau under the assumption that $c(\sfw)\ne s(\sfw)$ because, then, $J(\sfw) = A_{{}_{\overline{\sfw},-\mu(\sfw)^{-1}}}$.

\begin{thm} 
\label{thm.Af.ASR}
Suppose $k$ is algebraically closed and $\fchar(k) \ne 2,3$.
Let $\l \in k$, $f \in S^3V$, and let $E \subseteq \PP^2$ denote the zero locus of $f$.
\begin{enumerate}
\item
$A_{\!{}_{f,\l}}$ is 3-Calabi-Yau when $E$ is smooth.  
  \item 
$A_{\!{}_{f,\l}}$ is 3-Calabi-Yau when $E$ is $\vert\!\vert\!\vert\,$, $\leftrightline \!\!\!\!\!\!\vert\!\vert\,\,$, $\varhexstar\,$, 
$ \vert\!  \largecircle$, or $\bigcurlywedge$.  
 \item{}
If $E$ is $\largetriangleup$, $\varnothing\,$, or $\propto$, 
then $A_{\!{}_{f,\l}}$ is 3-Calabi-Yau except for two values of $\l$. 
\end{enumerate}
%Case (2) occurs if and only if $f \ne 0$ and $H^2(f)=0$. Case (3) occurs if and only if $H^2(f)\ne 0$ and the zero locus of $f$ is not smooth.
\end{thm}
\begin{pf}
%%% Whether the Hessian of a function in $SV$ is zero or not does not depend on the choice of basis for $V$.  Therefore, 
% to prove 
% the last two sentences in the statement of 
% the theorem 
%%% it suffices to prove 
% them 
%%% for the $f$'s in Table \ref{table.symm.relns.for.Jw}. This is easily done.
If $f=0$, then $A_{\!{}_{f,\l}}$ is the polynomial ring on 3 indeterminates so is 3-Calabi-Yau.
Hence, for the rest of the proof we assume that $f \ne 0$. 

Fix a $3 \times 1$ matrix $\sfx$ whose entries are a basis for $V$. Let $\sM$ be the unique $3 \times 3$ matrix such that
$\sfw_0-\l\widehat{f}=\sfx^\sT\sM\sfx$.  The entries in $\sM$ belong to $V$; we write $\overline{\sM}$ for $\sM$ when 
we treat its entries as elements of $SV$.

By Corollary \ref{cor.thm.3.2}, $A_{\!{}_{f,\l}}$ is standard in the sense of \S\ref{ssect.standard}.
Thus, by Theorem \ref{thm.ATV1.1}, $A_{\!{}_{f,\l}}$ is Artin-Schelter regular of dimension 3 if and only if   
$\rank(\overline{\sM}) \ge 2$ at all points of 
$\PP^2$.

If $\theta \in \GL(V)$, then $A_{\!{}_{f,\l}} \cong A_{\!{}_{\theta(f),\l\det(\theta)^{-1}}}$ by Proposition \ref{prop.Af.basis.change}. 
It therefore suffices to prove the theorem for  $A_{\!{}_{\theta(f),\l}}$ for some $\theta \in \GL(V)$.  This allows us to assume that $f$
is one of the polynomials in Table \ref{table.symm.relns.for.Jw}. 

(1)
Suppose $E$ is a smooth cubic. Replacing $f$ by a suitable $\theta(f)$, $\theta \in \GL(V)$, we can assume that $f= xyz+\frac{1}{3}\b(x^3+y^3+z^3)$ for some $\b \in k$ such that $\b^3\notin \{0,-1\}$. A straightforward calculation shows that $A_{\!{}_{f,\l}}$ has defining relations
\begin{align*}
ayz+bzy+cx^2 \; = \; & 0, \\
azx+bxz+cy^2  \; = \; & 0, \\
axy+byx+cz^2  \; = \; & 0, 
\end{align*}
where $a=\frac{1}{2}(\l-2)$, $b=\frac{1}{2}(\l+2)$, and $c=\b\l$. These are the relations at \cite[(1.4)]{ATV1}. As remarked at
 \cite[p.38]{ATV1}, this algebra fails to be Artin-Schelter regular  if and only if $a^3=b^3=c^3$ or two of $a,b,c$ are zero.

Since $\b \ne 0$, it is impossible for two of $a,b,c$ to be zero. Thus, $A_{\!{}_{f,\l}}$ is Artin-Schelter regular
of dimension 3 unless $(\l-2)^3=(\l+2)^3=(2\l \b)^3$.  If $(\l,\b)$ is a solution to these equations, then $(\l^2, \b^3)=(-\frac{4}{3},-1)$. Since $\b^3 \ne -1$ we conclude that $A_{\!{}_{f,\l}}$ is  Artin-Schelter regular of dimension 3 and therefore 3-Calabi-Yau.

(2)
Changing basis, we can, and do, assume that $f$ is one of the first five polynomials in Table \ref{table.symm.relns.for.Jw}.  
Arguing as in Lemma \ref{lem.22.entry-1}, it suffices to  show that the common zero locus of the $2 \times 2$ minors 
\begin{equation}
\label{eq.minors}
\begin{cases}
\phantom{x}  \hbox{$\frac{1}{4}$} \l^2 (f_{yy}f_{zz}-f_{yz}^2)+x^2 & \\
\phantom{x} \hbox{$\frac{1}{4}$} \l^2  (f_{zz}f_{xx}-f_{zx}^2)+y^2  &\\
\phantom{x} \hbox{$\frac{1}{4}$} \l^2  (f_{xx}f_{yy}-f_{xy}^2)+z^2 &
\end{cases}
\end{equation}
of $\overline{\sM}$ is empty. 

Suppose these three minors vanish at $p=(a,b,c) \in \PP^2$.
Since $f$ is either $x^3$, $x^2y$, $x^3+y^3$, $x^2z+xy^2$, or $x^2z+\frac{1}{3}y^3$, we have
$f_{zz}(p)=f_{yz}(p)=0$. From the first equation in (\ref{eq.minors}), we see that $a=0$.  Therefore $f_{zz}(p)=f_{zx}(p)=0$; from the 
second equation, we see that $b=0$. It follows that $f_{yy}(p)=f_{xy}(p)=0$; therefore $c=0$. We 
conclude that the common zero locus of these three minors is empty. Hence $A_{\!{}_{f,\l}}$ is  
 Artin-Schelter regular of dimension 3.

(3) 
In this case we can, and do, assume that $f=xyz+g$ where $g$ is either $0$, $\frac{1}{3} x^3$, or $\frac{1}{3} (x^3+y^3)$. 
By Lemma \ref{lem.22.entry-1},  $A_{\!{}_{f,\l}}$ is  3-Calabi-Yau except for two values of $\l$.
\end{pf}

\begin{cor}
\label{cor.Af.ASR}
Let $f \in S^3V$ and $\l \in k$. 
\begin{enumerate}
 \item 
If  $A_{\!{}_{f,\l}}$ is not 3-Calabi-Yau, then it is isomorphic to $k\langle x,y,z\rangle$ modulo the relations 
\begin{enumerate}
  \item 
  $xy=yz=zx=0$, or
  \item
 $xy=yz-x^2=zx=0$, or  
  \item 
  $xy=yz-x^2=zx-y^2=0$.  
\end{enumerate} 
  \item 
  $A_{\!{}_{f,\l}}$ is 3-Calabi-Yau if and only if it is a domain.
\end{enumerate} 
\end{cor}
\begin{pf}
(1)
Suppose $A_{\!{}_{f,\l}}$ is not 3-Calabi-Yau. 
Then we are in case (3) of Theorem \ref{thm.Af.ASR} and after a change of variables we can, and do, assume that $f=xyz+g$ where
$g \in \{0, \, \hbox{$\frac{1}{3}$} x^3, \, \hbox{$\frac{1}{3}$} (x^3+y^3)\}$. By Lemma \ref{lem.22.entry-1}(1) and its proof, $\l^2=4$, 
and $A_{\!{}_{f,\l}}$ is defined by the relations in (\ref{22.entry.relns}).

Suppose $\l=-2$. Then $(\a,\b)=(2,0)$ so the relations in (\ref{22.entry.relns}) become
$$
2 yz\, -   \, g_{xx}x \; = \;  2 zx\, -   \, g_{yy}y    \; = \;  2 xy  \; = \; 0.
$$
When $g=0$, the relations are $yz=zx=xy=0$; i.e., those in (a). If $g=\frac{1}{3}x^3$,  the relations    are 
$yz-x^2=zx=xy=0$; i.e., those in (b). When $g=\frac{1}{3}(x^3+y^3)$,  the relations  are $yz-x^2=zx-y^2=xy=0$; i.e., those in (c). 

Suppose $\l=2$. Then $(\a,\b)=(0,2)$ so the relations in (\ref{22.entry.relns}) are
$$
-2 zy\, -   \, g_{xx}x \; = \;  -2 xz\, -   \, g_{yy}y    \; = \;  -2 yx  \; = \; 0.
$$
When $g=0$ the relations  are $zy=xz=yx=0$. When $g=\frac{1}{3}x^3$,  the relations are $zy+x^2=xz=yx=0$. When
$g=\frac{1}{3}(x^3+y^3)$,  the relations are $zy+x^2=xz+y^2=yx=0$. 

{\it However,} the algebras obtained when $\l=2$ are isomorphic to those obtained when $\l=-2$. More precisely, $A_{\!{}_{f,2}} \cong A_{\!{}_{f,-2}}$: 
when $g$ is 0 or $\frac{1}{3}(x^3+y^3)$,  the linear map $\theta:V \to V$ defined by $\theta(z)=-z$,  
$\theta(x)=y$, and $\theta(y)=x$, extends to an isomorphism  $A_{\!{}_{f,2}} \to A_{\!{}_{f,-2}}$; 
when $g=\frac{1}{3}x^3$,  $\theta(z)=y$,  $\theta(x)=\sqrt {-1}x$, and $\theta(y)=z$, 
extends to an isomorphism  $A_{\!{}_{f,2}} \to A_{\!{}_{f,-2}}$.

(2)
Part (1) of this corollary shows that if $A_{\!{}_{f,\l}}$ is not 3-Calabi-Yau, then it is not a domain.
\end{pf}

\section{The classification when $c(\sfw) = s(\sfw)$}
\label{ssect.cw=sw}
In \S\ref{ssect.cw=sw} we assume $c(\sfw) = s(\sfw)$.

\subsection{}
If $\overline{\sfw}=0$ or, equivalently, if $E=\PP^2$, then $J(\sfw) \cong TV$  which is not Calabi-Yau because $\dim_k(V) \ne 1$. 
For the rest of \S\ref{ssect.cw=sw}, we assume that $\overline{\sfw} \ne 0$; i.e., $E$ is a cubic divisor on $\PP^2$.
 
\subsection{}
Since $J(\sfw)=J(c(\sfw))$ we can, and do, assume that $\sfw=c(\sfw)$.  

\subsection{}
Because $\sfw=s(\sfw)$,  $\sfw = \widehat{\overline{\sfw}}= \widehat{f}$.  
By Lemma \ref{lem.pd.hat}, the algebra $J(\sfw)$ in the next result is $J(\widehat{f}\,)$ and is therefore standard if the kernel of the map 
$V^{\otimes 2} \to J(\sfw)$ has dimension 3. 
A direct computation using the data in Table 2 shows that the kernel has dimension $<3$ if and only if  the zero locus of 
$f$ is $\PP^2$ or $\vert\!\vert\!\vert$ or $\leftrightline\!\!\!\!\!\vert\!\vert\;$ or $\varhexstar$ if and only if $H(\overline{\sfw})=0$.

\begin{prop}
\label{prop.not.CY}
Let $\sfw \in V^{\otimes 3}$ be such that $c(\sfw)=s(\sfw)$. 
Suppose that $\overline{\sfw} \ne 0$ and let $E \subseteq \PP^2$ be the subscheme $\{\overline{\sfw}=0\}$. 
\begin{enumerate}
\item
If $E$ is $\vert\!\vert\!\vert\,$, $\leftrightline  \!\!\!\!\!\vert\!\vert\,$, $\varhexstar$, $\vert\!  \largecircle$, or $\bigcurlywedge\,$, then $J(\sfw)$ is not a domain so is not 3-Calabi-Yau.
\item
If $E$ is $\largetriangleup$, $\varnothing\,$, or $\propto$, then $J(\sfw)$ is 3-Calabi-Yau and therefore a domain.
\item
If $E$ is smooth, then
\begin{enumerate}
  \item 
$J(\sfw)$ is a domain if and only if $j(E)\ne 0$; 
\newline
if  $j(E)\ne 0$, then $J(\sfw)$ is 3-Calabi-Yau;
  \item 
if $j(E)=0$, then $J(\sfw) \cong k\langle x,y,z\rangle/(x^2,y^2,z^2)$ and is not 3-Calabi-Yau;
\item
$J(\sfw)$ is 3-Calabi-Yau if and only if $j(E) \ne 0$.
\end{enumerate} 
 \item
If $J(\sfw)$ is not  3-Calabi-Yau, then it is not a domain.   
\end{enumerate}
\end{prop}
\begin{pf}
By Proposition \ref{prop.Jw.relns}, the hypothesis $c(\sfw) = s(\sfw)$ implies that $J(\sfw)$ is defined by the relations
$
 \widehat{\{y,z\}}_{\overline{\sfw}} \;= \;  \widehat{\{z,x\}}_{\overline{\sfw}} \;= \;    \widehat{\{x,y\}}_{\overline{\sfw}}  \;= \; 0.
$
Thus, if we write $f=\overline{\sfw}$, 
\begin{equation}
\label{c=s.algebras}
J(\sfw) \; \cong \; \frac{k \langle x,y,z\rangle}{(\widehat{f_x},\,\widehat{f_y}, \, \widehat{f_z})}
\end{equation}
by Lemma \ref{lem.pd.hat}. (The linear span of $f_x$, $f_y$, and $f_z$, does not depend on the choice of basis for $V$.)

As explained in \S\ref{ssect.cubics}, we can, and do, assume that $f$ is one of the polynomials in Table \ref{table.symm.relns.for.Jw}.

(1)
In this case, $J(\sfw)$ is one of the first five algebras in Table \ref{table.symm.relns.for.Jw}. One sees immediately from a glance at Table \ref{table.symm.relns.for.Jw} that these algebras are not domains. 
 
 (2)
Suppose $E$ is  $\largetriangleup$, $\varnothing$, or $\propto$. 
Let $f$ respectively be $xyz$, $xyz+\frac{1}{3}x^3$, or $xyz+\frac{1}{3}(x^3+y^3)$, and let
$\sM$ respectively be the matrix 
\begin{equation}
\label{M.for.symm.relns}
\begin{pmatrix} 0 & z & y \\ z & 0 & x \\ y & x & 0 \end{pmatrix}, \qquad
\begin{pmatrix} x & z & y \\ z & 0 & x \\ y & x & 0 \end{pmatrix}, \qquad
\begin{pmatrix} x & z & y \\ z & y & x \\ y & x & 0 \end{pmatrix}.
\end{equation}
Let $\sfx=(x,y,z)^\sT$. In each case the entries in $\sM\sfx$ are a basis for the space of relations defining $J(\sfw)$. 
For each $\sM$ it is easy to see that $\sfx^\sT\sM=(\sM\sfx)^\sT$ and the common zero locus of the $2 \times 2$ minors of 
$\sM$ is empty. Thus, by Theorem \ref{thm.ATV1.1},  $J(\sfw)$ is Artin-Schelter regular of dimension 3 and  therefore 3-Calabi-Yau
by Corollary \ref{cor.3ASR=3CY}.

(3)
Suppose $E$ is smooth. By changing the basis for $V$ we can, and do, assume that 
$
f= x^3+y^3+z^3 + 3 \l xyz
$
where $\l^3 \ne -1$. The defining relations for $J(\sfw)$ are 
\begin{equation}
\label{eq.some.sklyanin.algs}
\begin{cases}
2x^2 + \l yz + \l zy \; = \; 0, &\\
2y^2 + \l zx + \l xz \; =  \; 0, & \\
2z^2 + \l xy + \l yx \; =  \; 0. &
\end{cases}
\end{equation}
By \cite[Prop. 2.3]{PP} for example, the isomorphism class of $E$ is determined by its $j$-invariant which is
$$
j(E) \; = \; \left (\frac{3\l(8-\l^3)}{1+\l^3}\right )^3.
$$ 

(3a)
Suppose $j(E) \ne 0$.  
By  \cite[p.38]{ATV1}, the algebra with relations (\ref{eq.some.sklyanin.algs}) is a 3-dimensional Artin-Schelter regular algebra and therefore a domain.

(3b)
Suppose $j(E)=0$.\footnote{This is, in some sense the most symmetric case: if $k=\CC$, then $E \cong \CC/\ZZ+\ZZ\xi$ where $\xi$ is a primitive $6^{\th}$ root of unity so the lattice is hexagonal.}
Hence the relations for $A$ are $x^2=y^2=z^2=0$ so $A$ is not a domain and therefore not a 3-dimensional Artin-Schelter regular algebra.

(3c)
This follows from (3a) and (3b).

(4) 
Suppose that $J(\sfw)$ is not 3-Calabi-Yau.  If $E$ is singular, then $E$ is $\vert\!\vert\!\vert\,$, $\leftrightline  \!\!\!\!\!\vert\!\vert\,$, $\varhexstar$, $\vert\!  \largecircle$, or $\bigcurlywedge\,$ by (2), so $J(\sfw)$ is not a domain by (1).  If $E$ is smooth, then $j(E)=0$ by (3c) so  $J(\sfw)$ is not a domain by (3a). 

%%% Let $\cA$ denote the set of $A$ of the form (\ref{c=s.algebras}) that are 3-Calabi-Yau.
%%% Let $\cB$ denote the set of $A$ of the form (\ref{c=s.algebras}) that are not domains. 
% Taken together, (1) and (2) show that if $E$ is singular, then $J(\sfw)$ is either 3-Calabi-Yau or not a domain. 
%%% in $\cA \cup \cB$. 
% If $E$ is smooth and $j(E) = 0$, then $J(\sfw)$ is not a domain 
%%% $A \in \cB$ 
% by (3b). If  $E$ is smooth and $j(E)\ne 0$, then $J(\sfw)$ is 3-Calabi-Yau 
%%% $A \in \cA$ 
% by (3a).  
% Thus, every $J(\sfw)$ of the form (\ref{c=s.algebras}) is either 3-Calabi-Yau or not a domain. 
%%% belongs to $\cA \cup \cB$. 
\end{pf}

\begin{theorem}
 \label{thm.CY}  
Let $\sfw \in V^{\otimes 3}-\{0\}$.
\begin{enumerate}
  \item
  $J(\sfw)$ is 3-Calabi-Yau if and only if it is a domain.   
  \item 
Up to isomorphism of algebras, there are exactly nine $J(\sfw)$'s that are not  3-Calabi-Yau, namely the three algebras in \S\ref{ssect.Afl.not.ASR}, and the algebras in the rows of Table \ref{table.symm.relns.for.Jw} labelled with a $*$.
\end{enumerate}
\end{theorem} 
\begin{pf}
(1)
($\Rightarrow$)
If $J(\sfw)$ is  3-Calabi-Yau it is a domain by Corollary \ref{cor.3ASR=3CY}. 

($\Leftarrow$)
Suppose $J(\sfw)$ is not 3-Calabi-Yau. We will show it is not a domain.

Suppose $c(\sfw)\neq s(\sfw)$. By Theorem \ref{thm.3.2 }, $J(\sfw)\cong A_{{}_{f,\l}}$ for some $(f,\l)$. By Corollary \ref{cor.Af.ASR},
$A_{{}_{f,\l}}$ is 3-Calabi-Yau  if and only if it is a domain, so  $J(\sfw)$ is not a domain. 

Suppose $c(\sfw)=s(\sfw)$.  
Since $J(\sfw)$ is not 3-Calabi-Yau it is not a domain by Proposition \ref{prop.not.CY}(4).

(2)
Suppose $J(\sfw)$ is not 3-Calabi-Yau and therefore not a domain.

If $c(\sfw)\neq s(\sfw)$, then $J(\sfw)\cong A_{{}_{f,\l}}$ for some $(f,\l)$ by Theorem \ref{thm.3.2 }, so is 
one of the algebras in Corollary \ref{cor.Af.ASR}(1) which are exactly the algebras in  \S\ref{ssect.Afl.not.ASR}.

Suppose $c(\sfw) = s(\sfw)$.  Then $J(\sfw)$ is isomorphic to one of the algebras in Table \ref{table.symm.relns.for.Jw}
that is not a domain. Taken together, the non-domains in Table \ref{table.symm.relns.for.Jw} and the algebras
in Corollary \ref{cor.Af.ASR}(1) give nine isomorphism classes of $J(\sfw)$'s that are not 3-Calabi-Yau. 
\end{pf}

 In summary, because 3-dimensional AS-regular algebras are noetherian by \cite {ATV1}, we have the following result: 

\begin{cor} For $\sfw\in V^{\otimes 3}-\{0\}$, the following are equivalent:
\begin{enumerate}
\item{} $J(\sfw)$ is 3-Calabi-Yau. 
\item{} $J(\sfw)$ is noetherian and 3-Calabi-Yau. 
\item{} $J(\sfw)$ is a domain.
\item{} 
If $x, y\in V-\{0\}$, then $xy\ne 0$ in $J(\sfw)$. 
\end{enumerate}
\end{cor} 

The implication (4) $\Rightarrow$ (1) follows from the list of nine algebras in Theorem \ref{thm.CY}.
 
\subsection{Clifford algebras} 
\label{ssect.Clifford}

Let $R=k[X,Y,Z]$ be the polynomial ring on three indeterminates of degree 2 and let $M$ be a $3 \times 3$ symmetric matrix whose entries belong to 
$R$.
%${\rm span}_k\{X,Y,Z\}$. 
%There are unique matrices $M_1,M_2,M_3 \in M_3(k)$ such that $M=M_1X+M_2Y+M_3Z$.
The {\sf Clifford algebra} over $R$ associated to $M$ is the $R$-algebra
$$
A(M) \; :=\; \frac{R\langle x_1,x_2,x_3\rangle}{I}
$$
where $I$ is the ideal generated by the elements $x_ix_j+x_jx_i -M_{ij}$.
%$$   x_ix_j+x_jx_i = (M_1)_{ij}X+(M_2)_{ij}Y+(M_3)_{ij}Z, \qquad 1 \le i,j \le 3.  $$

Let $K=k(X,Y,Z)=\Fract(R)$. Since $K \otimes_R A(M)$ is the Clifford algebra over $K$ associated to $M$ it
has $K$-basis $\{1,x_1,x_2,x_3,x_1x_2,x_2x_3,x_3x_1,x_1x_2x_3\}$. 
It follows that these 8 elements are linearly independent over $R$. In $A(M)$, $x_jx_i=M_{ij}-x_ix_j$ and $x_i^2 \in R$,  so 
these 8 elements also generate $A(M)$ as an $R$-submodule. Hence $A(M)$ is a free $R$-module of rank 8 with the same basis.
%$\{1,x_1,x_2,x_3,x_1x_2,x_2x_3,x_3x_1,x_1x_2x_3\}$.  

The following example is relevant. If $a,b,c\in k$ and 
\begin{equation}
\label{symm.matrix.sklyanin}
M=\begin{pmatrix} 2X & cZ & bY \\ cZ & 2Y & aX \\ bY & aX & 2Z \end{pmatrix}, 
\end{equation}
then $A(M)=k[X,Y,Z,x_1,x_2,x_3]$ modulo the relations  
\begin{align*}
& \!\!\! \!  [X,-]   \; \; = \;  \; [Y,-]  \;  \; = \; \;  [Z,-] \;   \; = \;  \;  0,
\\
 \quad x_1^2 & = X, \quad \; \; x_2^2 =Y, \; \; \quad x_3^2=Z,
 \\
x_1x_2+x_2x_1=cZ, & \; \; \quad x_2x_3+x_3x_2 = aX, \;  \; \quad x_3x_1+x_1x_3 =bY.
\end{align*}
It follows that $A(M)$ is generated by $x_1,x_2,x_3$ modulo the relations
\begin{equation}
\label{cliff.relns.sklyanin1}
x_1x_2+x_2x_1 = cx_3^2,  \; \; \quad x_2x_3+x_3x_2 = ax_1^2, \;  \; \quad x_3x_1+x_1x_3 =bx_2^2, \quad [x_i^2,x_j]=0.
\end{equation}
%The next proof shows that the relations $[x_i^2,x_j]=0$ follow from the first three relations in  (\ref{cliff.relns.sklyanin1}).
%It then follows that $A(M)$ is generated by $x_1,x_2,x_3$ modulo the first three relations in (\ref{cliff.relns.sklyanin1}). 

There is some intentional repetition in the statement of the next result.

\begin{thm}
\label{thm.central.squares}
Let $\sfw \in V^{\otimes 3}-\{0\}$ be such that $c(\sfw)=s(\sfw)$. 
If $J(\sfw)$ is  3-Calabi-Yau, then there is a basis $\{x,y,z\}$ for $V$ such that 
\begin{enumerate}
  \item 
  $k[x^2,y^2,z^2]$ is a 3-dimensional polynomial ring,
  \item 
   $k[x^2,y^2,z^2]$ is in the center of $J(\sfw)$, and 
  \item 
  $J(\sfw)$ is a Clifford algebra of  rank 8 over  $k[x^2,y^2,z^2]$.  
  \item{}
    $J(\sfw)$ is a free $k[x^2,y^2,z^2]$-module of rank 8. 
\end{enumerate} 
\end{thm}
\begin{pf}
The hypothesis that $c(\sfw)=s(\sfw)$ implies that $J(\sfw)$ is one of the algebras in Table \ref{table.symm.relns.for.Jw}. 
The hypothesis that $J(\sfw)$ is 3-Calabi-Yau implies that $J(\sfw)$ is a domain and therefore isomorphic to one of the algebras in a row of 
Table \ref{table.symm.relns.for.Jw} that does not contain a $*$. Thus, up to a choice of basis $J(\sfw)$ is $k \langle x,y,z\rangle$
modulo the relations
\begin{equation}
\label{cliff.relns.sklyanin}
xy+yx=cz^2,  \; \; \quad yz+zy= ax^2, \;  \; \quad zx+xz =by^2
\end{equation}
for some $a, b, c\in k$. We can read off from Table \ref{table.symm.relns.for.Jw} that $abc \ne -1$.

It follows from (\ref{cliff.relns.sklyanin}) that
$$
 yx^2=(cz^2-xy)x=cz^2x-x(cz^2-xy) 
$$
whence $[y,x^2]=c[z^2,x]$. By symmetry, $[x,z^2]=b[y^2,z]$ and $[z,y^2]=a[x^2,y]$.
It now follows that $[y,x^2]=-abc[y,x^2]$. Since $abc \ne -1$, $[y,x^2]=0$. By symmetry, $[z,y^2]=[x,z^2]=0$ also.

Similarly, it follows from (\ref{cliff.relns.sklyanin}) that $y^2x=y(cz^2-xy)=cyz^2-(xy-cz^2)y$ whence $[y^2,x]=c[y,z^2]$. 
By symmetry, $[z^2,y]=a[z,x^2]$ and $[x^2,z]=b[x,y^2]$.
It now follows that $[y^2,x]=-abc[y^2,x]$. Since $abc \ne -1$, $[y^2,x]=0$. By symmetry, $[z^2,y]=[x^2,z]=0$ also.

Thus, %the relations (\ref{cliff.relns.sklyanin}) imply that 
$k[x^2,y^2,z^2]$ is contained in the center of $J(\sfw)$.

Therefore $J(\sfw)$ is generated by $x,y,z$ subject to the relations in (\ref{cliff.relns.sklyanin}) and $[x^2,-]=[y^2,-]=[z^2,-]=0$. 
Hence $J(\sfw)$ is isomorphic to the algebra $A(M)$ defined  just prior to the statement of this theorem. 
%On the other hand, if
%\begin{equation}
%\label{symm.matrix.sklyanin}
%M=\begin{pmatrix} 2X & cZ & bY \\ cZ & 2Y & aX \\ bY & aX & 2Z 
%\end{pmatrix},
%\end{equation} 
%then $A(M)$ is generated by $x_1,x_2,x_3$ modulo the six relations
%\begin{equation}
%\label{cliff.relns.sklyanin2}
%x_1x_2+x_2x_1=cx_3^2,  \; \; \quad x_2x_3+x_3x_2 = ax_1^2, \;  \; \quad x_3x_1+x_1x_3 =bx_2^2, \quad [x_i^2,x_j]=0. 
%\end{equation}
%The above proof shows that if $abc \ne -1$, then the relations $[x_i^2,x_j]=0$ follow from the first three relations. Thus, $J(\sfw)$ is the Clifford algebra over the polynomial ring $k[X,Y,Z]$ associated to the matrix in (\ref{symm.matrix.sklyanin}).   
It follows from this isomorphism that $k[x^2,y^2,z^2]$ is a polynomial ring on 3 variables.
\end{pf}

\section{The classification when $c(\sfw) \ne s(\sfw)$}
\label{sect.cw.not-=.sw}

\subsection{}
This section shows that the right-most column of  Table \ref{main.table} is correct. 
We therefore assume $c(\sfw) \ne s(\sfw)$ in \S\ref{sect.cw.not-=.sw}. By Theorem \ref{thm.3.2 }, 
$J(\sfw) \; \cong \; A_{{}_{\overline{\sfw},-\mu(\sfw)^{-1}}}$.
If $\overline{\sfw}=0$, then $J(\sfw)$ is the commutative polynomial ring on 3 indeterminates so is 3-Calabi-Yau. 
Thus, for the remainder of  \S\ref{sect.cw.not-=.sw} we assume that $\overline{\sfw} \ne 0$. 
This implies that $E$, the zero locus of $\overline{\sfw}$, is a cubic divisor.

\subsubsection{}
\label{S.12}
\label{S.32}
If  $E$ is smooth, then $J(\sfw)$ is  3-Calabi-Yau   by Theorem \ref{thm.Af.ASR}(1). 

\subsubsection{}
When $E$
is $\vert\!\vert\!\vert\,$, $\leftrightline  \!\!\!\!\!\vert\!\vert\,$, $\varhexstar\,$, $\vert\!  \largecircle$, or $\bigcurlywedge\,$,  
then $J(\sfw)$ is 3-Calabi-Yau by Theorem \ref{thm.Af.ASR}(2).

\subsubsection{}
We will show that the entries in Table \ref{main.table} for $\largetriangleup$, $\varnothing\,$, and $\propto$, are correct after the next lemma.

\begin{lem}
\label{lem.change.vars}
Let $x_1,\ldots,x_n$ and $X_1,\ldots,X_n$ be ordered bases for $V$, and $f \in SV$. 
Let $\nabla_{\!{}_x}^2(f)$ and $\nabla_{\!{}_X}^2(f)$ be the Hessian matrices for $f$ with respect to $x_1,\ldots,x_n$ and $X_1,\ldots,X_n$ respectively. Let $H_{\!{}_x}(f)$ and $H_{\!{}_X}(f)$ be the determinants of  $\nabla_{\!{}_x}^2(f)$ and $\nabla_{\!{}_X}^2(f)$, respectively. Let $A \in \GL(n)$ be the unique matrix such that
$$
(X_1,\ldots,X_n)=(x_1,\ldots,x_n)A^\sT.
$$
\begin{enumerate}
  \item 
  $\nabla_{\!{}_x}^2(f)  \; = \; A^{\sT}\nabla_{\!{}_X}^2(f) A$. 
  \item 
  $H_{\!{}_x}(f)=(\det A)^2 H_{\!{}_X}(f)$.
  \item 
  If $\theta \in \GL(V)$ and $H(\theta(f))$ and $H(f)$ are computed with respect to the same basis for $V$, 
  then $H(\theta(f))=(\det \theta)^2 \, \theta\big(H(f)\big)$. 
\end{enumerate}
\end{lem}
\begin{pf}
(1)
Let $A=(a_{ij})$. Since $X_i=a_{i1}x_1+\cdots+a_{in}x_n$, $a_{ij}=\pd X_i/\pd x_j$. Therefore
$$
\frac{\pd f}{\pd x_j} \;= \;   \sum_{p=1}^n \frac{\pd X_p}{\pd x_j} \frac{\pd f}{\pd X_p}  \;= \;  
 \sum_{p=1}^n a_{pj} \frac{\pd f}{\pd X_p}. 
 $$
 The $ij^{\th}$ entry in $\nabla_{\!{}_x}^2(f)$ is 
 $$
 \frac{\pd^2 f}{\pd x_i \pd x_j} \;= \;   \sum_{p=1}^n a_{pj} \frac{\pd }{\pd x_i}\Bigg( \frac{\pd f}{\pd X_p} \Bigg) \;= \;  
 \sum_{p=1}^n a_{pj}  \sum_{q=1}^n a_{qi}    \frac{\pd^2 f}{\pd X_q X_p}
 $$
 which is the $ij^{\th}$-entry in $A^{\sT}\nabla_{\!{}_X}^2(f) A$. The result follows.
 
 (2)
It follows from (1) that  the determinant of  $\nabla_{\!{}_x}^2(f)$ is $\det(A)^2$ times the determinant of $\nabla_{\!{}_X}^2(f)$.

(3)
Define $X_i:=\theta(x_i)$ and let $A \in \GL(n)$ be the unique matrix such that
$(X_1,\ldots,X_n)=(x_1,\ldots,x_n)A^\sT$. By (1), 
$$
\nabla_{\!{}_x}^2(\theta(f))  \; = \; A^{\sT}\nabla_{\!{}_X}^2(\theta(f)) A \; = \; A^{\sT}\theta\big(\nabla_{\!{}_x}^2(f) \big) A.
$$
Since $\det(\theta)=\det(A)$, the result is obtained by taking the determinant of both sides of the previous equality and applying (2).
\end{pf}

\begin{prop}
\label{prop.22.entry-1}
Let $\{x,y,z\}$ be a basis for $V$ and let $\theta \in \GL(V)$.
Let $f  \in S^3V$,   and let $E \subseteq \PP^2$ be the zero locus of $f$.
\begin{enumerate}
  \item 
  If $E=\largetriangleup$, then $A_{\!{}_{f,\l}}$ is 3-Calabi-Yau if and only if $\l^2H(f) \ne 8f$.
  \item 
If $E$ is $\varnothing$ or  $\propto$, then $A_{\!{}_{f,\l}}$ is 3-Calabi-Yau if and only if 
the zero locus of $\l^2H(f) +24 f$ is not a triangle.
\end{enumerate} 
\end{prop}
\begin{pf} 
If $E=\largetriangleup\,$, let $g=0$. If $E=\varnothing\,$, let $g=\frac{1}{3}x^3$. If $E=\propto$, let $g=\frac{1}{3}(x^3+y^3)$.
By the discussion in \S\ref{ssect.cubics}, there is $\theta \in \GL(V)$ such that $f=\theta(xyz+g)$. 
Since $H(xyz+g)=2(xyz-3g)$,  
$$
H(f)=H(\theta(xyz+g)) =2(\det \theta)^2\theta (xyz-3g).
 $$
By Proposition \ref{prop.Af.basis.change}, $A_{\!{}_{f,\l}}  = A_{\!{}_{\theta(xyz+g),\l}} \cong A_{\!{}_{xyz+g,\l\det(\theta)}}$. 
Therefore $A_{\!{}_{f,\l}}$ is 3-Calabi-Yau if and only if $A_{\!{}_{xyz+g,\l\det(\theta)}}$ is. Thus, by 
Lemma \ref{lem.22.entry-1}(2), $A_{\!{}_{f,\l}}$ is 3-CY if and only if $ \l^2\det(\theta)^2 \ne 4$. The following equivalences are obvious: 
\begin{align*}
 \l^2\det(\theta)^2 \ne 4 
&  \; \Longleftrightarrow \; \l^2(\det \theta)^2  \, \theta\big(2xyz-6g\big) \ne 4  \, \theta\big(2xyz-6g\big) 
\\
&  \; \Longleftrightarrow \; \l^2H (f) \ne 8  \, \theta\big(xyz-3g\big) 
\\
& \; \Longleftrightarrow \;  \l^2H(f) \ne 8f -32 \, \theta(g),  \quad \qquad \;\, (*)
\\
& \; \Longleftrightarrow \;  \l^2H(f) +24 f  \ne 32 \, \theta(xyz). \qquad (**)
\end{align*}

(1)
If $E=\largetriangleup\,$, $g=0$ so, by $(*)$, $A_{\!{}_{f,\l}}$ is  3-Calabi-Yau if and only if   
$\l^2H(f) \ne 8f$.  

(2)
We now assume $E$ is $\varnothing$ or  $\propto$. By $(**)$, (2) holds if the following conditions on 
 $\l$ are equivalent:
\begin{enumerate}
  \item[(a)] 
$\l^2H(f) +24 f  \ne 32 \, \theta(xyz)$
  \item[(b)]
  the zero locus of $\l^2H(f) +24 f$ is not a triangle. 
\end{enumerate}  
Certainly, (b) implies (a). 

Suppose (a) is true. To show (b) holds it suffices to show there is only one triangle in the pencil generated by 
$H(f)$ and $f$. Let $g$ be $\frac{1}{3}x^3$ or $\frac{1}{3}(x^3+y^3)$. 
It is easy to see that the only triangle in the pencil generated by $xyz$ and $g$ is $xyz=0$. Since $H(xyz+g)=2(xyz-3g)$ it follows that there is only one triangle in the pencil generated by $H(xyz+g)$ and $xyz+g$. Hence there is only one triangle in the pencil generated by $\theta(xyz+g)=f$ and $\theta(H(xyz+g))= (\det\theta)^{-2}H( \theta(xyz+g))=(\det\theta)^{-2}H(f)$; i.e.,
there is only one triangle in the pencil generated by $f$ and $H(f)$.  
\end{pf}

\begin{cor}
\label{cor.22.entry}
Let $\sfw \in V^{\otimes 3}-\{0\}$ and suppose that $c(\sfw)\ne s(\sfw)$. 
 Let $E$ be the zero locus of $\overline{\sfw}$. 
\begin{enumerate}
  \item 
If $E=\largetriangleup$, 
then  $J(\sfw)$ is  3-Calabi-Yau if and only if $H(\overline{\sfw}) \ne 8\mu(\sfw)^2 \overline{\sfw}$. 
  \item 
If $E$ is $\varnothing$ or $\propto$, then $A_{\!{}_{f,\l}}$ is 3-Calabi-Yau if and only if 
the zero locus of $H(\overline{\sfw}) +24\mu(\sfw)^2 \overline{\sfw}$ is not a triangle.
\end{enumerate}
In both (1) and (2) there are exactly two values of $\mu(\sfw)$ for which $J(\sfw)$ is not 3-Calabi-Yau.
\end{cor}
\begin{pf}
Since  $c(\sfw)\ne s(\sfw)$,  $J(\sfw) = A_{\!{}_{\overline{\sfw},-\mu(\sfw)^{-1}}}$. 
So $J(\sfw)$ is 3-CY if and only if the condition in Proposition \ref{prop.22.entry-1} for 
$A_{\!{}_{f,\l}}$ to be  3-CY holds with $f$ replaced by $\overline{\sfw}$ and $\l$ replaced by $-\mu(\sfw)^{-1}$.  
\end{pf}

\section{The point scheme} 
 \label{sect.pt.var}

\subsection{}
A 3-dimensional Artin-Schelter regular algebra $A$ determines, and is determined by, a 
pair $(D, \theta)$ where $D$ is either $\PP^2$ or a cubic divisor in $\PP^2$ and $\theta$ is an automorphism of $D$.   It is common to call $D$ the {\it point scheme} of $A$ because it is a fine moduli space for certain $A$-modules called point modules---see \cite{ATV1} for more information. 
If $J(\sfw)$ is 3-Calabi-Yau and $D \ne \PP^2$, Theorem \ref{thm.pt.sch} determines $D$.

\begin{lem}
\label{lem.det.sym+skew-sym}
Let $A$ and $B$ be $3 \times 3$ matrices with entries in a commutative $k$-algebra $R$. Let $A_{ij}$ and $B_{ij}$  denote their $ij^{\th}$ entries, etc.
Let $\l \in k$. If $A$ is symmetric and $B$ skew-symmetric, then $\det(\l A+B) = ( \det A)\l^3+c \l$ where
\begin{align*}
c    & \; = \;  B_{23}^2A_{11}   +  B_{31}B_{23} A_{12} +  B_{12}B_{23}A_{13}
\\
& \qquad +B_{31}B_{23}A_{12} +B_{31}^2A_{22} +B_{12}B_{31}A_{23}
\\
& \qquad \qquad
 +B_{12}B_{23}A_{13} +B_{12}B_{31} A_{23} +B_{12}^2A_{33}.
\end{align*}
\end{lem}
\begin{pf}
Let $a,b,c,d \in R$ be such that $\det(\l A+B) =a\l^3+b\l^2+c\l+d$.  Now
\begin{align*}
\det(\l A+B) & \; = \; \det(\l A^\sT+B^\sT)
\\
& \; = \; \det(\l A-B)
\\
&  \; = \; -\det(-\l A+B) \qquad \hbox{because the matrix has odd size}
\\
&  \; = \; a\l^3 - b\l^2+c\l -d.
\end{align*}
Therefore $b=d=0$. It is obvious that $a=\det(A)$. 

Write $A_i$ and $B_i$ for the $i^{\th}$ rows of $A$ and $B$, and $A_{ij}$ and $B_{ij}$ for the $ij^{\th}$ 
entries of $A$ and $B$, respectively. 
The coefficient of $\l$ is 
\begin{align*}
& A_1 \wedge B_2 \wedge B_3  \; + \; B_1 \wedge A_2 \wedge B_3 \; + \;  B_1 \wedge B_2 \wedge A_3 
\\
& \qquad  \;= \;  
\left| \begin{matrix} A_{11} & A_{12} & A_{13} \\ -B_{12} & 0 &  B_{23} \\ B_{31} & -B_{23} & 0 \end{matrix} \right|
\; + \; 
\left| \begin{matrix} 0 & B_{12} &  -B_{31} \\ A_{12} & A_{22} & A_{23} \\  B_{31} & -B_{23} & 0 \end{matrix} \right|
\; + \; 
\left| \begin{matrix} 0 & B_{12} &  -B_{31} \\  -B_{12} & 0 & B_{23}  \\ A_{13} & A_{23} & A_{33}  \end{matrix} \right|
\\
& \qquad  \;= \; A_{11}  B_{23}B_{23}  + A_{12} B_{31}B_{23} +A_{13}B_{12}B_{23} 
\\
& \qquad \qquad +A_{12}B_{31}B_{23} +A_{22}B_{31}B_{31} +A_{23}B_{12}B_{31}
\\
& \qquad \qquad  \qquad +A_{13}B_{12}B_{23} +A_{23}B_{12}B_{31} +A_{33}B_{12}B_{12}.
\end{align*}
The proof is complete. 
\end{pf}

\begin{thm} 
\label{thm.pt.sch}
Let $\sfw \in V^{\otimes 3}$. If $J(\sfw)$ is 3-Calabi-Yau, then its point scheme is the subscheme 
$H(\overline{\sfw})\, + \, 24 \, \mu(\sfw)^2 \overline{\sfw}\; =\; 0$ of $\PP^2$.
\end{thm}
\begin{pf}
Fix a basis $\{x,y,z\}$ for $V$ and define $\sfx:=(x,y,z)^\sT$.

Suppose $c(\sfw)=s(\sfw)$. Then
$$
J(\sfw) \cong \frac{k\langle x,y,z\rangle}{(\widehat{f_x}, \widehat{f_y},\widehat{f_z})}
$$
where $f=\overline{\sfw}$. Since $f_x \in S^2V$, Lemma \ref{lem.g.hat} tells us that $\widehat{f_x}= f_{xx}x+f_{xy}y+f_{xz}z$.
There are similar expressions for $\widehat{f_y}$ and $\widehat{f_z}$. It follows that the relations for $J(\sfw)$ are the entries in 
$\sM\sfx$ where   
$$
\sM \; = \; 
\begin{pmatrix}
 f_{xx} & f_{xy} & f_{xz}   \\
 f_{yx} & f_{yy} & f_{yz}   \\
  f_{zx} & f_{zy} & f_{zz}   
\end{pmatrix}
\; = \; \nabla^2(f).
$$
By \cite[Cor. 3.13]{ATV1},  the point scheme for $J(\sfw)$ is the zero locus of $\det(\overline{\sM})$; i.e.,  the zero locus of 
$H(f)=H(\overline{\sfw})$. Since $\mu(\sfw)=0$ when $c(\sfw)=s(\sfw)$, the theorem is true when $c(\sfw)=s(\sfw)$.

Suppose  $c(\sfw) \ne s(\sfw)$. Then $J(\sfw) \cong A_{\!{}_{f,\l}}$ where $f=\overline{\sfw}$ and  $\l=-\mu(\sfw)^{-1}$. 
By \S4.1, the defining relations of   $A_{\!{}_{f,\l}}$ are the entries in $\sM\sfx$ where 
$$
\sM=\begin{pmatrix} 
\frac{\l}{2}f_{xx} & \frac{\l}{2}f_{xy}+z & \frac{\l}{2}f_{xz}-y \\
\frac{\l}{2}f_{yx}-z & \frac{\l}{2}f_{yy} & \frac{\l}{2}f_{yz}+x \\
\frac{\l}{2}f_{zx}+y & \frac{\l}{2}f_{zy}-x & \frac{\l}{2}f_{zz} \\
\end{pmatrix}.
$$
The point scheme for $J(\sfw)$ is given by $\det(\overline{\sM})=0$.
We can write $\sM=\l A+B$ where $A$ is the symmetric matrix $\frac{1}{2}\nabla^2 (\overline{\sfw})$ and  
$$
B=\begin{pmatrix}  0 & z &-y \\ -z & 0 &  x \\  y & -x & 0 \end{pmatrix}.
$$
By Lemma \ref{lem.det.sym+skew-sym}, $\det(\overline{\sM}) = \frac{1}{8}\l^3 H(\overline{\sfw}) + \frac{1}{2}\l c$ where 
\begin{align*}
c & \;=\;  x^2 f_{xx}  +xy f_{xy}  +xzf_{xz} +yxf_{yx}+y^2f_{yy} +yzf_{yz}+zxf_{zx}+zyf_{zy} +z^2f_{zz}
\\
 & \;=\; x(x f_{xx}  +y f_{xy}  +zf_{xz})\,  +\, y(xf_{yx}+yf_{yy} +zf_{yz}) \, + \, z(xf_{zx}+yf_{zy} +zf_{zz})
 \\
 & \;=\; 2xf_x\,  +\, 2yf_y \, + \, 2zf_z
  \\
 & \;=\; 6f.
\end{align*}
Thus, the point scheme for $J(\sfw)$ is the locus $ \frac{1}{8}\l^3 H(\overline{\sfw}) + 3\l \overline{\sfw}=0$. 
\end{pf}
 
To end, we describe the 3-CY algebras for which $c(\sfw)=s(\sfw)$ and $E$ is an elliptic curve.
The members of that class of algebras are classified in the proof of 
Proposition \ref{prop.not.CY}(3). Up to isomorphism, they are $k\langle x,y,z\rangle$ modulo relations of the form (\ref{eq.some.sklyanin.algs}) and, in (\ref{eq.some.sklyanin.algs}), $\l^3 \notin \{0,8\}$ because the $j$-invariant of $E$ is non-zero. Thus, the  3-CY algebras $J(\sfw)$ for which $c(\sfw)=s(\sfw)$ and $E$ is an elliptic curve are, up to isomorphism, precisely the algebras $A$ described in the next result.

 \begin{prop}
Let $c \in k$ and suppose that $c^3 \notin\{0,1\}$. Let $A$ be $k\langle x,y,z\rangle$ modulo the relations
\begin{align*}
cx^2+ yz+ zy \; = \; & 0 \\
cy^2+ zx+ xz \; = \; & 0 \\
cz^2+ xy+ yx \; = \; & 0. 
\end{align*}
The point variety for $A$ is given by $c(x^3+y^3+z^3) - (2+c^3)xyz =0$ and its automorphism is translation by a point of order 2.
\end{prop}
\begin{pf}
Let $g = c(x^3+y^3+z^3) - (2+c^3)xyz$ and write $D$ for the zero locus of $g$. 
As remarked at \cite[(1.6)]{ATV1}, $D$ is the point scheme for $A$ and the canonical automorphism, $\theta$, of $D$ is given by 
\begin{align*}
\theta(x,y,z) \; = \; & (xz-cy^2,\, yz-c x^2, \, c^2 xy-z^2)
\\
 \; = \; & (x^2-c^2yz,\,c z^2-xy, \,cy^2-xz).
\end{align*}
We use two formulas because there are points on $D$ where one of the formulas gives $\theta(x,y,z)=(0,0,0)$. 
As at \cite[p.38]{ATV1}, we fix the group law $(D,+,o)$ with identity $o=(1,-1,0)$. Straightforward calculations show that $\theta(o)=(1,1,c)$ and $\theta^2=\id_D$.  

We will show that $\theta$ is translation by $p:=(1,1,c)$.

Both $g$ and its Hessian vanish at $o$, so $o$ is an inflection point, and therefore 
$$
\hbox{three points $p,q,r \in D$ are collinear} \; \Longleftrightarrow \; p+q+r=o.
$$

Let $(\a,\b,\c) \in D$.  
Since $\c x + \c y - (\a+\b)z$ vanishes at $(\a,\b,\c)$, $(\b,\a,\c)$, and $o$, 
$$
-(\a,\b,\c)=(\b,\a,\c).
$$

Since $g_x(p)=g_y(p)$, $o$ is on the line $g_x(p)x+g_y(p)y+g_z(p)z=0$. This is the tangent line to $D$ at $p$ so $p+p=0$;
i.e., $p$ is 2-torsion.   

The claim is that if $(\a,\b,\c) \in D$, then $\theta(\a,\b,\c)=(\a,\b,\c) +(1,1,c)$. Since $(1,1,c)$ is 2-torsion and $-(\a,\b,\c)=(\b,\a,\c)$, 
this is equivalent to the claim that
\begin{equation}
\label{transl.autom}
\theta(\a,\b,\c)+(\b,\a,\c) +(1,1,c) =o.
\end{equation}
Both $(\b,\a,\c)$ and $(1,1,c)$ lie on the line $(\a c - \c)x+(\c-c\b)y+(\b-\a)z=0$. 
A calculation shows that $(\a\c-c\b^2,\b\c-c\a^2,c^2\a\b-\c^2)$ also lies on this line; but the latter point is 
$\theta(\a,\b,\c)$. Hence (\ref{transl.autom}) holds.
\end{pf}

 \begin{prop}
Let $\sfw \in V^{\otimes 3}$. Suppose that $c(\sfw)=s(\sfw)$ and that the zero locus of 
$\overline{\sfw}$ is an elliptic curve whose $j$-invariant is non-zero. 
Then the point variety for $J(\sfw)$ is a smooth elliptic curve and its automorphism is translation by a point of order 2.
\end{prop}

 \end{document}